%% file: main.tex
\def\year{2019}\relax
\documentclass[letterpaper]{article} 
\usepackage{aaai19}  
\usepackage{times}  
\usepackage{helvet}  
\usepackage{courier}  
\usepackage{url}  
\usepackage{graphicx}  
\usepackage[noend]{algpseudocode}
\usepackage{algorithm,amsmath}
\usepackage{subfig}
\usepackage{mathtools}
\usepackage{newtxtext,newtxmath}

\usepackage{amsthm}
\usepackage{comment}
\usepackage[utf8]{inputenc}
\usepackage[dvipsnames]{xcolor}
\newcommand*{\R}{\mathbb{R}}
\newcommand*{\Po}{\text{Prox}}

\newcommand*{\E}{\mathbb{E}}

\newcommand{\norm}[1]{\left\lVert#1\right\rVert}
\newcommand{\Iprod}[2]{\left\langle #1,#2\right\rangle}
\newcommand\myeq[2]{\mathrel{\stackrel{{{#1}}}{#2}}}

\newcommand{\Initialize}{\textbf{Initialize:}{\,}}
\newcommand{\Input}{\textbf{Input:}{\,}}
\newcommand{\Output}{\textbf{Output:}{\,}}

\newtheorem{theorem}{Theorem}
\newtheorem{lemma}{Lemma}

\newtheorem{definition}{Definition}

\newtheorem{remark}{Remark}
\newtheorem{assumption}{Assumption}
\newcommand{\keepcomment}{1}
\AtBeginDocument{\ifnum\keepcomment=1
  \excludecomment{comment}
\else
  \includecomment{comment}
\fi}
\begin{document}
%
\title{Asynchronous Delay-Aware Accelerated Proximal Coordinate Descent for Nonconvex Nonsmooth Problems}
\author{
Ehsan Kazemi, Liqiang Wang\\
Department of Computer Science,
University of Central Florida\\
}
\maketitle
\begin{abstract}
Nonconvex and nonsmooth problems have recently attracted considerable attention in machine learning. However, developing efficient methods for the nonconvex and nonsmooth optimization problems with certain performance guarantee remains a challenge. Proximal coordinate descent (PCD) has been widely used for solving optimization problems, but the knowledge of PCD methods in the
nonconvex setting is very limited. On the other hand, the asynchronous proximal coordinate descent (APCD) recently have received much attention in order to solve large-scale problems. However, the accelerated variants of APCD algorithms are rarely studied. In this paper, we extend APCD method to the accelerated algorithm (AAPCD) for nonsmooth and nonconvex problems that satisfies the sufficient descent property, by comparing between the function values at proximal update and a linear extrapolated point using a delay-aware momentum value. To the best of our knowledge, we are the first to provide stochastic and deterministic accelerated extension of APCD algorithms for general nonconvex and nonsmooth problems ensuring that for both bounded delays and unbounded delays every limit point is a critical point. By leveraging Kurdyka-Łojasiewicz property, we will show linear and sublinear convergence rates for the deterministic AAPCD with bounded delays. 
Numerical results demonstrate the practical efficiency of our algorithm in speed. 
\end{abstract}

\noindent 

\section{Introduction}
For many machine learning and data mining applications, efficiently solving the optimization problem with nonsmooth regularization is important. In this paper, we focus on the following composite optimization
problem of machine learning model with nonsmooth regularization term as
\begin{equation}\label{problem}
    \min_{x\in\R^m}\,F(x) = f(x) + g(x)
\end{equation}
where $f:\R^m \to\R$ captures the empire risk which is smooth and possibly nonconvex, and $g:\R^m \to\R$, corresponding to the regularization term, reduces to a finite-sum
\begin{equation}
g(x) = \sum_{j=1}^m g_j(x_j)
\end{equation}
where each $g_j$ can be nonconvex. 

Many problems on \eqref{problem} correspond to convex model that can be efficiently 
optimized by first order algorithm, in particular accelerated proximal gradient (APG) methods which is proven to be efficient for the class of convex functions. However, many real applications require the problems to be nonconvex. The nonconvexity might originate either from function $f(x)$ or the regularization function. This type of problems is popular in machine learning, for example, sparse logistic regression \cite{liu2009large}, and sparse multi-class classification \cite{blondel2013block}. On the other hand regarding the nonsmooth regularization terms, proximal gradient methods often address solving optimization problems with nonsmoothness. The proximal operator is defined as following
\[
\Po_{\eta g_j}(y) =  \text{arg}\min_{x\in\R^m}\frac{1}{2\eta}\norm{x-y}^2+g_j(x_j)
\]
where $\eta > 0$, and $\norm{\cdot}$ is $l_2$-norm. If the proximal operator does not have an analytic solution, an algorithm should be used to  solve the proximal operator which might be inexact. In this paper we consider only algorithms which use exact proximal mapping.

While the new algorithms for problem \eqref{problem} provide both good theoretical convergence and empirical performances, the investigations on them
were mainly conducted in the sequential setting.
 In the current big data era, we need to design algorithms to deal with
very large scale problems  ($m$ is large). In this case, we need 
to eliminate sequential updates which usually take too much costly idle time. This necessitates parallel computation which will not use synchronization to wait for all others and share their updates. Recently asynchronous parallelization have received huge
successes due to its potential to vastly speed up algorithms \cite{dean2012large,recht2011hogwild}. We design and analyze an asynchronous parallel implementations of the accelerated proximal coordinate descent algorithms with bounded and unbounded delays for nonconvex nonsmooth problems, which is not well studied in the literature, to the best of our knowledge.
\subsection{Contributions}
The main contributions of this paper are
summarized as follows. We first propose the basic stochastic and deterministic variants of asynchronous accelerated  proximal
 coordinate descent algorithm for nonconvex problems.   
By construction of Lyapunov functions, we show that the limit points of the sequences generated by AAPCD are critical points of the problem \eqref{problem} for both bounded delays and unbounded delays. This is one of the first convergence results for a method with acceleration which alleviates the bottleneck of unbounded delays for nonsmooth nonconvex functions. The convergence studies for AAPCD, through a novel perspective, characterize the stepsize based on the momentum parameter. This fills the void in previous analyses such as \cite{li2015accelerated,yao2017efficient}, where the effect of the exact value of the momentum parameter on the acceleration of convergence were not observed. As the stability of the algorithm is highly affected by asynchronism, by allowing negative momentum for high staleness values we will show the reduction in the objective function will be increased significantly and accelerates convergence. In particular, we characterize the momentum parameter in the sense that increasing the stepsize would involve decreasing of the momentum parameter, while it will provide comparable asymptotic convergence in terms of the violation of first-order optimality
conditions. We will show that by requiring momentum, a fixed stepsize could be chosen for unbounded delays.
   
By leveraging different cases of Kurdyka-Łojasiewicz property of the objective function, we establish the linear and sub-linear convergence rates of the function value sequence generated by the deterministic AAPCD with bounded deterministic delays and they match the synchronous results. In all the cases investigated in this paper, the independence assumption between blocks and delays is avoided.

We provide numerical experiments to demonstrate  the performance of our stochastic AAPCD algorithm on various
large-scale real-world datasets. The results outperforms other asynchronous stochastic algorithms reported in literature such as ASCD \cite{liu2015asynchronous} and
AASCD \cite{fang2018accelerating}. It also shows that
AAPCD can achieve good speedup on large-scale
real-world datasets and provide significantly faster convergence to a reasonable accuracy than
competing options, while still providing favorable asymptotic accuracy. 

\section{Related Works}
{\bf{Proximal Gradient Algorithms:}} Proximal gradient methods for nonsmooth regularization are among the most important methods for solving composite optimization problems. There have been accelerated exact proximal gradient variants.
Specifically, for convex problems, the authors in \cite{beck2009fast} displayed basic accelerated  proximal gradient (APG) method which extends
Nesterov's accelerated methods for solving single smooth convex function \cite{nesterov1983method}. They proved that APG displays the non-asymptotic convergence rate $O(\frac{1}{k^2})$, where $k$ is the number of iterations. 

For extensions to nonconvex settings, \cite{ghadimi2016accelerated} studied the condition that only the regularization term could be nonconvex,
and proved the convergence rate of APG method. \cite{boct2016inertial} established the convergence of proximal method when 
 $f(x)$ and $g(x)$ could be nonconvex. \cite{li2015accelerated} focused on first-order algorithms and by exploiting KL property they proved that APG algorithm can converge to a stationary point in different rates. Recently, in \cite{gu2016inexact} and \cite{li2017convergence} several accelerated proximal methods were studied, and sublinear and linear rates under different cases of the KL property for nonconvex problems were provided.

In addition to the above proximal gradient
methods, several stochastic optimization methods were developed for solving composite problems see, e.g.,  proximal stochastic coordinate descent prox-SCD \cite{shalev2011stochastic}, prox-SVRG \cite{xiao2014proximal}, prox-SAGA \cite{defazio2014saga}, prox-SDCA \cite{shalev2014accelerated}. Under the assumption that the regularization term is block separable,  \cite{richtarik2014iteration} developed a randomized block-coordinate descent method. An accelerated variant of this method is studied in \cite{lin2015accelerated}. All these stochastic methods require convexity of $f$, or even stronger assumptions. 

For nonconvex problems, \cite{ghadimi2016accelerated} generalized an accelerated SGD method to solve nonconvex but smooth minimization problems. Stochastic  variance reduction methods for nonconvex problems were investigated in \cite{allen2016variance,reddi2016stochastic}.  Furthermore, proximal variance reduction methods for general nonconvex, nonsmooth problems are proposed in \cite{reddi2016proximal,allen2017natasha}. Then, \cite{xu2015block} proposed a block stochastic gradient method for nonconvex and nonsmooth problems. 

{\bf Asynchronous Coordinate Descent:}
The asynchronous computation is much more
efficient than the synchronous computation. More recently, asynchronous parallel methods have been successfully applied to accelerate many optimization algorithms including stochastic  coordinate descent \cite{liu2015asynchronous}. We briefly review the works which are closely related to
ours as follows. ASCD can provide linear and sublinear convergence rates \cite{liu2015asynchronous,avron2015revisiting}. Similar results were established for asynchronous SGD \cite{recht2011hogwild}, and stochastic variance reduction algorithms \cite{reddi2015variance,leblond2017asaga}. A study of ASCD for unbounded delays has been performed in \cite{sun2017asynchronous}, however the results are restricted only to Lipschitz differentiable functions. Some asynchronous algorithms particularly outperform conventional ones.
In \cite{meng2016asynchronous}, authors integrated momentum acceleration and variance reduction techniques to accelerate asynchronous SGD. Several accelerated schemes for asynchronous coordinate descent and SVRG using momentum compensation techniques were proposed in \cite{fang2018accelerating}. Recently, \cite{hannah2018texttt} analyzed an asynchronous accelerated block coordinate descent algorithm with optimal complexity which converges linearly to a solution for strongly convex functions.

However, to the best of our knowledge, there is no study on the asynchronous parallel versions of accelerated proximal coordinate descent algorithms for nonconvex nonsmooth objective functions. 

\section{Preliminaries and Assumptions}
We describe our asynchronous accelerated proximal coordinate descent for nonconvex problems in Algorithm \ref{AAPCD-Algo}. Compared to the regular proximal coordinate descent step, AAPCD takes an extra linear extrapolation step depending on the value of the current ages of $\hat{y}^k$, which is called also delay and denoted by $d_k$.  In order to  compute the delay $d_k$, we use a scalar counter to denote the weights at iteration $k$, starting from $k=0$, and with each update we increment the counter by one. We allow each worker to record the iteration $i$ when reading the weights and we let $k$ to denote the iteration when the same worker updating the weights. Then the actual delay $d_k$ is $d_k = k - i$. If delay is greater than the threshold $T_1$, we consider adding negative momentum to extrapolate a new iterate. We further show that adding such a
momentum for large delays have the effect of decreasing Lyapunov function over iterations. For acceleration, AAPCD only accepts the new extrapolated iterate when the objective function value is sufficiently decreased. It is important to note that the threshold $T_1$ can adaptively change during the iterations. From practical point of view there is a need to know how to select the parameter $T_1$. We will address this question later when we present the analyses of convergence.
It will be shown that accumulation points of sequences generated by AAPCD will converge to stationary points of $F$. 
In the step 5 of Algorithm \ref{AAPCD-Algo}, at iteration $k$, the block gradient $\nabla_{j_k} f$ is computed at the delay iterate $\hat{y}^k$, which is assumed to be some earlier state of $y^k$ in the shared memory with the delay $d_k$. The delay iterate $\hat{y}^k$ can be formulated as 
\begin{equation}
\hat{y}^k = y^k - \sum_{h\in I(k)} (y^{h+1}-y^h)
\end{equation}  
where $I(k) \in \{k-1,\ldots,k-d_k\}$ is a subset of previous iterations. From the proximal update for AAPCD, we have $x^{k+1}_j=y^{k}_j$ for $j\neq j_k$. We also assume $\beta = \max_{k}\{\beta_k\geq 0\}$ and $\beta' = \max_{k}\{\beta_k< 0\}$. We let $\Gamma_{k}^{r}$ be the set of iterations from $k$ to $r$ with $d_k> T_1$  and $y^{k+1} = v^k$, ${\Gamma^c}_{k}^{r}$ denote the set of iterations from $k$ to $r$ with $d_k\leq T_1$ and $y^{k+1} = v^k$, and ${\Gamma^0}_{k}^{r}$ denote the set of iterations from $k$ to $r$ with $y^{k+1} = x^{k+1}$.
\begin{algorithm}
\caption{Asynchronous Accelerated Proximal Coordinate Decent (AAPCD)}
\label{AAPCD-Algo}
\begin{algorithmic}[1]
\State\Input The stepsize $\eta$, threshold $T_1$
\State\Initialize $y^0 \in \R^m$
\For{ $k=0, 1,\ldots, R$ }
\State Randomly choose $j_k$ from
$\{1,\ldots,m\}$
\State $x_{j_k}^{k+1}=\Po_{j_k,\eta g_{j_k}}{(y^k-\eta\nabla_{j_k} f(\hat{y}^{k}}))$ and $x_{j}^{k+1} = y_{j}^{k}$ for $j\neq j_k$
\If {$d_k \leq T_1$}{\text{ choose} $\beta_k > 0$}
\State $v^k_{j_k}=x^{k+1}_{j_k}+\beta_k(x^{k+1}_{j_k}-y^{k}_{j_k})$ and $v_{j}^{k} = y_{j}^{k}$ for $j\neq j_k$
\Else {\text{ choose} $\beta_k < 0$}
\State $v^k_{j_k}=x^{k+1}_{j_k}+\beta_k(x^{k+1}_{j_k}-y^{k}_{j_k})$ and $v_{j}^{k} = y_{j}^{k}$ for $j\neq j_k$
\EndIf
\If {$F(x^{k+1})\leq F(v^k)$}
\State $y_{j_k}^{k+1}=x_{j_k}^{k+1}$
\Else 
\State{$y^{k+1}_{j_k}=v^k_{j_k}$} 
\EndIf
 \EndFor
 \State\Output ${y}_{R+1}$
\end{algorithmic}
\end{algorithm}

By studying different cases of KL property we will show that AAPCD will decrease the function value properly at the initial point. For the deterministic AAPCD with deterministic bounded staleness, we prove the linear and sublinear convergence rate by exploiting different cases of KL property. 

In the following we first introduce some tools for analyzing asynchronous algorithms, and then describe the
assumptions on the problem \eqref{problem} that we assume in this
paper.

For analysis of the stochastic algorithm, we let $\mathcal{F}_k$ denote the sigma algebra generated by $\{y^0,\ldots, y^k\}$ . We denote the total expectation by $\E$ and the expectation over the stochastic variable $d_k$ by $\E_{d_k}$. Function $g(x)$ is lower semicontinuous at point $x_0$ if $\lim\inf_{x\to x_0} g(x) \geq g(x_0)$. Throughout this paper, we assume each $g_j$ in problem \eqref{problem} is lower semicontinuous. A point $x\in \R^m$ is said a critical point of function $F$ if $0\in\partial F(x)$. 
The following Uniformized KL property is a powerful tool to analyze the first order descent algorithms.
 \begin{definition}[Uniformized KL Property]
A function $f:\R^m\to (-\infty,\infty]$ is said to satisfy the Uniformized KL property if for every compact set $\Omega\subset \text{dom}~\partial f$ on which $f$ is constant,  there exists $\epsilon, \gamma\in (0,+\infty]$ and $\phi \in \Phi_{\gamma}$, such that for all $\hat{u}\in\Omega$ and all $u\in \{u\in\R^m:\text{dist}_{\Omega}(u) <\epsilon\}\cap\{u\in\R^m: f(\hat{u})<f(u)<f(\hat{u})+\gamma\}$, the following inequality holds 
\[
\phi'(f(u)-f(\hat{u}))\text{dist}_{\partial f(u)}(0) \geq 1
\]
where  $\Phi_{\gamma}$ stands for a class of function $\phi:[0,\gamma)\to\R^+$ satisfying: (1) $\phi$ is concave and $C^1$ on $(0,\gamma)$; (2) $\phi$ is continuous at $0$, $\phi(0)=0$; and (3) $\phi'(x)>0$, for all $x\in (0,\gamma)$.
\end{definition}
By \cite[Lemma 6]{bolte2014proximal}, if function $f$ is lower semicontinuous and satisfies KL property at every point of $\Omega$, then it satisfies the Uniformized KL property. 
All semi-algebraic functions satisfy the KL property. Specially, the desingularising function $\phi(t)$ of semi-algebraic functions can be chosen to take the form $\phi(t) = \frac{C}{\theta}t^\theta$ with
$\theta\in (0,1]$. In particular, typical semi-algebraic functions include real polynomial functions, $\norm{x}_p$ with $p \geq 0$, rank, etc.

We make the following assumptions on the problem \eqref{problem} in this paper.
\begin{assumption}\label{P2-assu1}
Function $f$ and each $g_j$ are proper and lower semicontinuous; $\inf_{x\in\mathbb{R}^m} F(x) > -\infty$; the sublevel set
$\{x \in \mathbb{R}^d : F (x) \leq \alpha\}$ is bounded for all $\alpha\in\mathbb{R}$.
\end{assumption}
\begin{assumption}\label{P2-assu2} Function $f$ is continuously differentiable and the gradient $\nabla f$ is $L$-Lipschitz continuous.
\end{assumption}
To prove the limit points of $\{y^k\}$ generated by AAPCD are stationary points, we need a new assumption:
\begin{assumption} \label{unbdd-determinstic-rule}
For AAPCD, it is assumed that there exists $K \in \mathbb{N}$ such that for all $k \in \mathbb{N}$,
we have $\{1, \ldots,m\} \subseteq \{j_{k+1},\ldots,j_{k+K}\}$.
\end{assumption}
The goal of our paper is to provide a comprehensive analysis for AAPCD for both bounded and unbounded delays to justify the overall advantages of AAPCD.

\section{AAPCD with Bounded Delays}
In this section we analyze the convergence of Algorithms \ref{AAPCD-Algo} for bounded delays, i.e., we assume $d_k\leq \tau$ for all $k$ and for a fixed number $\tau$. 
Define the Lyapunov function $G$ as
\[
G(x^k):=G(x^k,y^k,\ldots,y^{k-\tau}) = F(x^k) + \xi_k
\]
where the sequence $\{\xi_k\}_{k\in\mathbb{N}}$, defined by
\[
\xi_k := \frac{L^2\tau}{2C}\sum_{h=k-\tau+1}^k(h-k+\tau)\norm{y^h-y^{h-1}}^2 
\]
with $C>0$ is a constant to be determined later.
In the lemma below, we present an inequality which states for a proper stepsize, AAPCD can provide sufficient descent in our Lyapunov function. 
\begin{lemma}\label{P2-lemma1} Suppose Assumption \ref{P2-assu2} hold. Given $\eta > 0$, we have
{
\begin{equation}\label{p2-lemma-eq-v1}
\begin{split}
\E[&G(x^{k+1})]\leq \E[G(y^{k})]\\
&~~-\left(\frac{1}{2\eta}-\frac{L}{2}-{L\tau(1+\beta_k)}\right)\E\norm{x^{k+1}-y^k}^2.\\
\end{split}
\end{equation}
}
\end{lemma}
We characterize the convergence of
AAPCD. Our first result characterizes the behavior of
the limit points of the sequence generated by AAPCD. Based on the lemma, we show that the sequence $\{y^k\}$ generated by AAPCD approaches critical points of the general nonconvex problem \eqref{problem}.
\begin{theorem}\label{P2-theo1}
Let Assumptions \ref{P2-assu1}-\ref{unbdd-determinstic-rule}  hold for the
problem \eqref{problem}. Then with stepsize $\eta < \frac{1}{L+2LT_1(1+\beta)}$, and the momentum $-1 < \beta_k < \frac{1}{L\tau}(\frac{1}{2\eta}-\frac{L}{2})-1$ the sequence $\{y^k\}$ generated by AAPCD satisfies
\begin{enumerate}
\item $\{y^k\}$ is an almost surely bounded sequence and $\E\norm{y^{k+1}-y^k}\to 0$.

\item The set of limit points of $\{y^k\}$ forms a compact set, on which function $F$ is a constant $F^*$ and the sequences $\{F(y^k)\}$ and $\{G(y^k)\}$ converge to $F^*$.

\item All the limit points of $\{y^k\}$ are critical points of $F$, and $\E[\text{dist}_{\partial F(y^{k})}(0)]={o}(\frac{1}{\sqrt{k}})$.
\end{enumerate}
\end{theorem}
\begin{remark}
The connectedness and compactness of the set 
$\Omega$ of the limit points of $\{y^k\}$ is implied from $E\norm{y^{k+1} - y^k}\to 0$. Theorem \ref{P2-theo1} also states that the objective function on $\Omega$ containing the critical points remains constant. 
\end{remark}
\begin{remark}
Equation \eqref{p2-lemma-eq-v1} shows that the selection of negative $\beta_k$ for substantial staleness values would increase Lyapunov function reduction over an iteration. In the light of the bounds for the momentum term $\beta_k$ in Theorem \ref{P2-theo1}, we could realize an estimation of an upper bound for the threshold $T_1$ in AAPCD algorithm. The staleness bound $T_1$ should be large enough to allow positive $\beta_k$. For example if $\beta=\frac{1}{2}$, then we should have, $\frac{1}{L\tau}(\frac{1}{2\eta}-\frac{L}{2})-1\geq \frac{1}{2}$. Thus, by choosing $\eta=\frac{1}{L+4LT_1(1+\beta)}$, we obtain $T_1\geq \frac{3\tau}{4(1+\beta)}=\frac{\tau}{2}$. 
\end{remark}

The compact set $\Omega$ satisfies the requirements of the Uniformized KL property, and hence can be utilized to show the decrease of function values, depending on a certain exponent $\theta$ defined below. 

\begin{theorem}
Let the conditions of Theorem \ref{P2-theo1} hold. Suppose that $F$ satisfies the Uniformized KL property with desingularising function $\phi$  of the form $\phi(t) = \frac{e}{\theta} t^\theta$.  Let $F(x) = F^*$ for all
of the limit points of $\{x^k\}$ in AAPCD, and denote  $r_k = F(x^k) - F^*$. Then the sequence $\{r_k \}$ for $k$ large enough satisfies
\begin{enumerate}\label{P2-theo2}
\item If $\theta = 1$, and $x_0$ is chosen such that $r_0< \frac{1}{b_1e^2}$, then $r_k$ reduces to zero in finite steps;
\item If $\theta= \frac{1}{2}$, then $\E[r_{k+1}]\leq \frac{b_1e^2}{1+b_1e^2} \E[r_{0}]$,
\end{enumerate}
where 
\[
 b_1 = \frac{2(\frac{1}{\eta}+L)^2(K+1)+2L^2T_1(1+\beta)+2L^2T(1+\beta'')}{\left(\frac{1}{2\eta}-\frac{L}{2}-{L\tau(1+\beta)}\right)}
 \] 
with $\beta''=\max\{\beta',0\}$.
\end{theorem}
\begin{remark}
As $\beta''\leq 0$, the contribution of the delays greater than $T_1$ in the factor $b_1$, i.e., $2L^2T(1+\beta'')$ decreases, which indicates acceleration is possible with negative momentum term.
\end{remark}
\section{AAPCD with Unbounded Delays}
In this section, we allow the delay $d_k$ to be an unbounded stochastic variable, and extremely large delays in our algorithm are permitted. Depending on some limitations on the distribution of $d_k$, we can still prove convergence. For unbounded delay analysis, one approach is to consider a new bound for the distribution of the  end-behavior of $d_k$ to decay sufficiently fast as the iterations progress.

We emulate this solution in the following. In particular, we define fixed parameters $p_j$ related to probabilities of the delay such that $p_j \geq P(j(k) = j)$, for all $k$,  and $c_{k}:=\sum_{t=1}^\infty t(t+k)p_{t+k}$ with  $\sum_{k=0}^{\infty} c_k <\infty$.  For instance, we note that if ${p_j}$ have  the probability distributions with decay bound $p_j={\mathcal{O}}(j^{-t})$, $t>4$, then $\sum_{k=0}^{\infty} c_k$ is finite.

We define a more involved Lyapunov function $G$ as
\begin{equation}
G(x^k) := G(x^k,y^k,\ldots,y^0)  = F(x^k) + \xi_k
\end{equation}
where to simplify the presentation, we define $\xi_k$ which encompasses all terms
\[
\xi_k :=\frac{L^2}{2C} \sum_{h=1}^{k} c_{k-h}\norm{y^h-y^{h-1}}^2.
\] 
where $\frac{1}{C} > 0$ is a contraction rate to be defined later.

\begin{lemma}\label{stc-unbdd-P2-lemma1} Under Assumption \ref{P2-assu1}, for any $\eta > 0$, we have 
\begin{equation}\label{stc-unbdd-p2-lemma-eq}
\begin{split}
\E[&G(x^{k+1})]\leq \E[G(y^{k})]\\
&-\left(\frac{1}{2\eta}-\frac{L}{2}-{L(1+\beta_k){\sqrt{c_0}}}\right)\E\norm{x^{k+1}-y^k}^2.\\
\end{split}
\end{equation}
\end{lemma}
Now we characterize the behavior of
the limit points of the sequence generated by AAPCD with unbounded delays. 
\begin{theorem}\label{stc-unbdd-P2-theo1}
Let Assumptions \ref{P2-assu1}-\ref{unbdd-determinstic-rule} hold for the
problem \eqref{problem}. Then with stepsize $\eta<\frac{1}{L+2L\sqrt{c_{T_1}}(1+\beta)}$ and momentum $-1 < \beta_k < \frac{1}{L\sqrt{c_0}}(\frac{1}{2\eta}-\frac{L}{2})-1$,  the sequence $\{y^k\}$ generated by AAPCD satisfies
\begin{enumerate}
\item $\{y^k\}$ is an almost surely bounded sequence and $\E[\xi_k]\to 0$.

\item The set of limit points of $\{y^k\}$ forms a compact set, on which the functions $F$ is a constant $F^*$ and $\{F(y^k)\}$ and $\{G(y^k)\}$ converge to $F^*$.

\item All the limit points of $\{y^k\}$ are critical points of $F$.
\end{enumerate}
\end{theorem}
\begin{remark}
Lemma \ref{stc-unbdd-P2-lemma1} shows that the selection of negative $\beta_k$ for delays greater than $T_1$ would decrease Lyapunov function substantially over an iteration. The  bounds for $\beta_k$ in Theorem \ref{stc-unbdd-P2-theo1} imply an estimation of a lower bound for $c_{T_1}$. For example if $\beta=\frac{1}{2}$, then, we should have $\frac{1}{L\tau}(\frac{1}{2\eta}-\frac{L}{2})-1\geq \frac{1}{2}$. Hence, by selecting $\eta=\frac{1}{L+4L\sqrt{c_{T_1}}(1+\beta)}$, we obtain $c_{T_1}\geq \frac{9c_0}{16(1+\beta)^2}=\frac{c_0}{4}$.    
\end{remark}
Now by applying the Uniformized KL property we show Algorithm  \ref{AAPCD-Algo} decreases the objective value below that of $F(x_0)$.
\begin{theorem}\label{stc-unbdd-P2-theo2}
Let the conditions of Theorem \ref{stc-unbdd-P2-theo1} hold and $F$ satisfies the Uniformized KL property and the desingularising function has the form of $\phi(t) = \frac{e}{\theta}t^{\theta}$ with $e > 0$.  We denote $r_k = F(y^k) - F^*$, where $F^*$ is the function value on the set of limit points of $\{y^k\}$. Then for $k$ large enough the sequence $\{r_k \}$  satisfies 
\begin{enumerate}
\item If $\theta = 1$, and $x_0$ is chosen such that $r_0< \frac{1}{b_1e^2}$ then $r_k$ reduces to zero in finite steps;
\item If $\theta= \frac{1}{2}$, then $\E[r_{k+1}]\leq \frac{b_1e^2}{1+b_1e^2} \E[r_{0}]$,
\end{enumerate}
where 
\[
b_1 = \frac{(\frac{2}{\eta^2}+4L^2)+4L^2+{4L^2c_0}(1+\beta)}{\frac{1}{2\eta}-\frac{L}{2}-{L(1+\beta){\sqrt{c_0}}}}.
\] 
\end{theorem}
\section{Deterministic AAPCD}
In this section, we consider deterministic unbounded delays. Specifically, deterministic AAPCD is presented in Algorithm \ref{Deterministic-AAPCD-Algo}. The stochastic and deterministic AAPCD differ only on how the current coordinates are selected at each iteration. For this purpose, we assume the delay variable $d_k$ is deterministic, which allow extremely large delays in our algorithm. We will prove that a subsequence of points $\{y^k\}$  generated by deterministic AAPCD converges to a stationary point. Using KL property we will see that if $x^0$ is not a stationary point, Algorithm \ref{Deterministic-AAPCD-Algo} decreases
the objective value below that of $F(x^0)$. We also prove the rate of convergence for the deterministic algorithm with deterministic bounded delay by exploiting KL property,  which is unavailable in the stochastic setting for the Lyapunov function.
\begin{algorithm}
\caption{Deterministic AAPCD}\label{Deterministic-AAPCD-Algo}\begin{algorithmic}[1]
\State\Input The stepsize $\eta$, threshold $T_1$
\State\Initialize $y^0 \in \R^m$
\For{ $k=0, 1, 2,\ldots, R$ }
\State Choose $j_k$ from
$\{1,\ldots,m\}$
\State {$x_{j_k}^{k+1}=\Po_{j_k,\eta g_{j_k}}{(y^k-\eta\nabla_{j_k} f(\hat{y}^{k}}))$ and $x_{j}^{k+1} = y_{j}^{k}$ for $j\neq j_k$}
\If {$d_k \leq T_1$} {choose $\beta_k>0$}
\State {$v^k_{j_k}=x^{k+1}_{j_k}+\beta_k(x^{k+1}_{j_k}-y^{k}_{j_k})$ and $v_{j}^{k} = y_{j}^{k}$ for $j\neq j_k$}
\Else { choose $\beta_k <0$}
\State {$v^k_{j_k}=x^{k+1}_{j_k}+\beta_k(x^{k+1}_{j_k}-y^{k}_{j_k})$ and $v_{j}^{k} = y_{j}^{k}$ for $j\neq j_k$}
\EndIf
\If {$F(x^{k+1})\leq F(v^k)$} 
\State {$y_{j_k}^{k+1}=x_{j_k}^{k+1}$}
\Else 
\State{$y^{k+1}_{j_k}=v^k_{j_k}$} 
\EndIf
 \EndFor
 \State\Output ${y}_{R+1}$
\end{algorithmic}
\end{algorithm}

As recommended in \cite{sun2017asynchronous}, we set a sequence $\{\epsilon_i\}_{i\geq 0}$ and define the Lyapunov function $G$ which encompasses all terms to control unbounded delays
\begin{equation}
G(x^k) := F(x^k) + \xi_k
\end{equation}
where to simplify the presentation, we define
\[
\xi_k :=\frac{L^2}{2C} \sum_{h=1}^{\infty} \delta_{k-h}\norm{y^h-y^{h-1}}^2
\] 
with $\delta_i = \sum_{j=i}^\infty \epsilon_j$ such that $\sum_{j=0}^\infty \delta_j < \infty$ and $C > 0$ to be determined later.

\begin{lemma}\label{unbdd-P2-lemma1} Let Assumption \ref{P2-assu2} hold. For any $\eta > 0$, we have
\begin{equation}\label{unbdd-p2-lemma-eq}
\begin{split}
G(&x^{k+1})\leq G(y^{k})\\
&~~-\left(\frac{1}{2\eta}-\frac{L}{2}-\sqrt{{\delta_0}{\mu_{d_k}}}L(1+\beta_k)\right)\norm{x^{k+1}-y^k}^2\\
\end{split}
\end{equation}
where $\mu_{d_k} = \sum_{h=0}^{d_k-1}\frac{1}{\epsilon_{h}}$.
\end{lemma}
For any $T \geq \lim \inf d_k$ which can be arbitrarily large, let $S_T$ be the subsequence of $\mathbb{N}$ where the current delay is less than $T$. We will show the
points $x^k$, $k \in S_T$, have convergence guarantees. 
The following theorem for unbounded deterministic delay is parallel to Theorem \ref{stc-unbdd-P2-theo1}.
 
\begin{theorem}\label{unbdd-P2-theo1}
Suppose that Assumptions \ref{P2-assu1}-\ref{unbdd-determinstic-rule} hold. Then with stepsize $\eta = \frac{c}{L+2 \sqrt{{\delta_0}{\mu_{T_1}}}L(1+\beta)}$ for $c\in (0,1)$, and momentum $-1 < \beta_k < {\frac{\sqrt{\mu_{T_1}}}{c\sqrt{\mu_{d_k}}}}(1+\beta)-1$, we have, 
\begin{enumerate}
\item $\{y^k\}$ is a bounded sequence and $\xi_k\to 0$.

\item The function $F$ is constant on the set of limit points of $\{y^k\}$ and the sequences $\{F(y^k)\}$ and $\{G(y^k)\}$ converge to it.

\item For any subsequence $S_T$ generated by the deterministic AAPCD, all the limit points of $\{y^k\}_{k\in S_T}$ are critical points of $F$.
\end{enumerate}
\end{theorem}
\begin{remark}
Lemma \ref{unbdd-P2-lemma1} shows that the use of momentum for delayed gradient might gain no performance and have negative effects. Hence, to compensate this issue, we allow the selection of negative $\beta_k$ for high staleness values to maximize the reduction of
the Lyapunov function over an iteration. By taking the bounds in Theorem \ref{unbdd-P2-theo1} for the momentum term $\beta_k$ in to consideration, we could present an upper bound estimate for the threshold $T_1$ in AAPCD. The delay bound $T_1$ should be large enough to allow positive $\beta_k$. For example if we choose $\beta=1$, then we should have, ${\frac{\sqrt{\mu_{T_1}}}{c\sqrt{\mu_{d_k}}}}(1+\beta)-1\geq 1$. Therefore, $T_1$ must be large enough such that ${\mu_{T_1}}\geq \frac{4c^2\mu_{d_k}}{(1+\beta)^2}=c^2\mu_{d_k}$, for all $k$. 
\end{remark}
It is important to note that although Theorem \ref{unbdd-P2-theo1} shows a fixed step size works for deterministic AAPCD, however, in return the upper bound for momentum is adaptive to the current delay.

In the following theorem, it turns out that a subsequence of Algorithm \ref{Deterministic-AAPCD-Algo} can decrease the function value at $x_0$, depending on the parameter $\theta$ defined below. 
\begin{theorem}\label{unbdd-deter-subseq-conver}
Let conditions of Theorem \ref{unbdd-P2-theo1} hold and that $F$ satisfies the Uniformized KL property and the desingularising function has the form $\phi(s)=\frac{e}{\theta}t^\theta$, where $\theta\in (0, 1]$ and $e > 0$.  Let $F(x) = F^*$ for all $x \in\Omega$ (the set
of limit points), and denote $r_k = F(y^k) - F^*$. Then the sequence $\{r_k \}_{k\in S_T}$ for $k$ large enough satisfies 
\begin{enumerate}
\item If $\theta = 1$, and $x_0$ is chosen such that $r_0< \frac{1}{b_1c^2}$ then $r_k$ reduces to zero in finite steps;
\item If $\theta\in [\frac{1}{2},1)$, then for $k$ large enough $r_k\leq \frac{b_1e^2}{1+b_1e^2} r_{0}$;
\item If $\theta\in (0,\frac{1}{2})$, then  
$r_{k}\leq \left(\frac{1}{b_2{({1-2\theta})}+r_{{0}}^{2\theta-1}}\right)^{\frac{1}{1-2\theta}}
$
\end{enumerate}
where 
 \[
 b_1=\frac{2(\frac{1}{\eta}+L)^2+3(1+\beta)^2L^2T_1+2(1+\beta'')^2L^2T}{(\frac{1}{c}-1)\frac{L}{2}}
 \] 
 with $\beta''=\max\{\beta',0\}$ and $b_2 = \min(\frac{1}{b_1e^2R},\frac{r_{0}^{2\theta-1} (R^{\frac{2\theta-1}{2\theta-2}}-1)}{1-2\theta})$  for a fixed number $R\in(1, \infty)$. 
\end{theorem}
For the deterministic AAPCD with deterministic bounded delay $T$, we define $\epsilon_i=0$ for $i > T$ and we let $\tilde{G}(x^k,y^k,\ldots,y^{k-T})$
 denote the corresponding Lyapunov function. In the following $\tilde{G}(x)$ refers to $\tilde{G}(x,x,\ldots,x)$. We let $\Omega$ denote the set of stationary points of $F$. Since $\xi_k\to 0$, by Theorem \ref{unbdd-P2-theo1}, $\tilde{G}$ is constant on $\Omega$. We can derive convergence rates for
\begin{equation}
r_k = \tilde{G}(y^k)-F^*.
\end{equation}
\begin{theorem}\label{bd-deter-corro}
Assume the conditions of Theorem \ref{unbdd-deter-subseq-conver}, but only $\tilde{G}$ satisfies the Uniformized KL property and the desingularising function has the form $\phi(s)=\frac{e}{\theta}t^\theta$, where $\theta\in (0, 1]$ and $e > 0$. Then if the delay is bounded by $T$, the sequence $\{r_k \}$ for $k$ large enough  satisfies
\begin{enumerate}
\item If $\theta = 1$, then $r_k$ reduces to zero in finite steps;
\item If $\theta\in [\frac{1}{2},1)$, then $r_{k}\leq \left(\frac{b_1e^2}{1+b_1e^2}\right)^{\lfloor\frac{k-k_1}{T+K}\rfloor}r_{k_1}$ for $k_1$ large enough;
\item If $\theta\in (0,\frac{1}{2})$, then  $r_{k}\leq \left(\frac{1}{{\lfloor\frac{k-k_0}{T+K}\rfloor}b_2{({1-2\theta})}+r_{{0}}^{2\theta-1}}\right)^{\frac{1}{1-2\theta}}$,
\end{enumerate}
where 
 \begin{equation}
 \begin{split}
 b_1=&\frac{3}{(\frac{1}{c}-1)\frac{L}{2}}\bigg((\frac{1}{\eta}+L)^2+(1+\beta)^2L^2T_1\\&+(1+\beta'')^2L^2T+2(1+\beta)^2{L^2\mu_{T}} \delta_{0}\bigg)
 \end{split}
\end{equation}
 with $\beta''=\max\{\beta',0\}$ and $b_2 = \min(\frac{1}{b_1e^2R},\frac{r_{0}^{2\theta-1} (R^{\frac{2\theta-1}{2\theta-2}}-1)}{1-2\theta})$  for a fixed number $R\in(1, \infty)$.
\end{theorem}
The convergence rates in Theorem \ref{bd-deter-corro} match the results from \cite{davis2016asynchronous}, but they need the independence assumption between blocks and delays. If $T=0$ we obtain a synchronous version of the accelerated coordinate descent, and hence Theorem \ref{bd-deter-corro} implies the same rates as given in \cite{li2015accelerated} for nonconvex functions. 
\begin{remark}
The characterization of the factor $b_1$ in Theorems \ref{unbdd-deter-subseq-conver} and \ref{bd-deter-corro} is noticeable in a particular way that the delays greater than $T_1$ contribute to this factor. Since $\beta'' \leq 0$, it shows that applying negative momentum for high delay values could efficiently decrease the value of $b_1$ which results in acceleration.
\end{remark}
\begin{remark}
The KL property of $F$ is not necessarily sufficient to ensure that the Lyapunov function $G$ satisfies the KL property. However, since $G-F$ is semi-algebraic and the class of
semi-algebraic functions is closed under addition, it shows
that $G$ is semi-algebraic, which implies that $G$ is a KL function.
\end{remark}
\section{Numerical Results}
\begin{figure}
\centering
\includegraphics[width=0.6\linewidth]{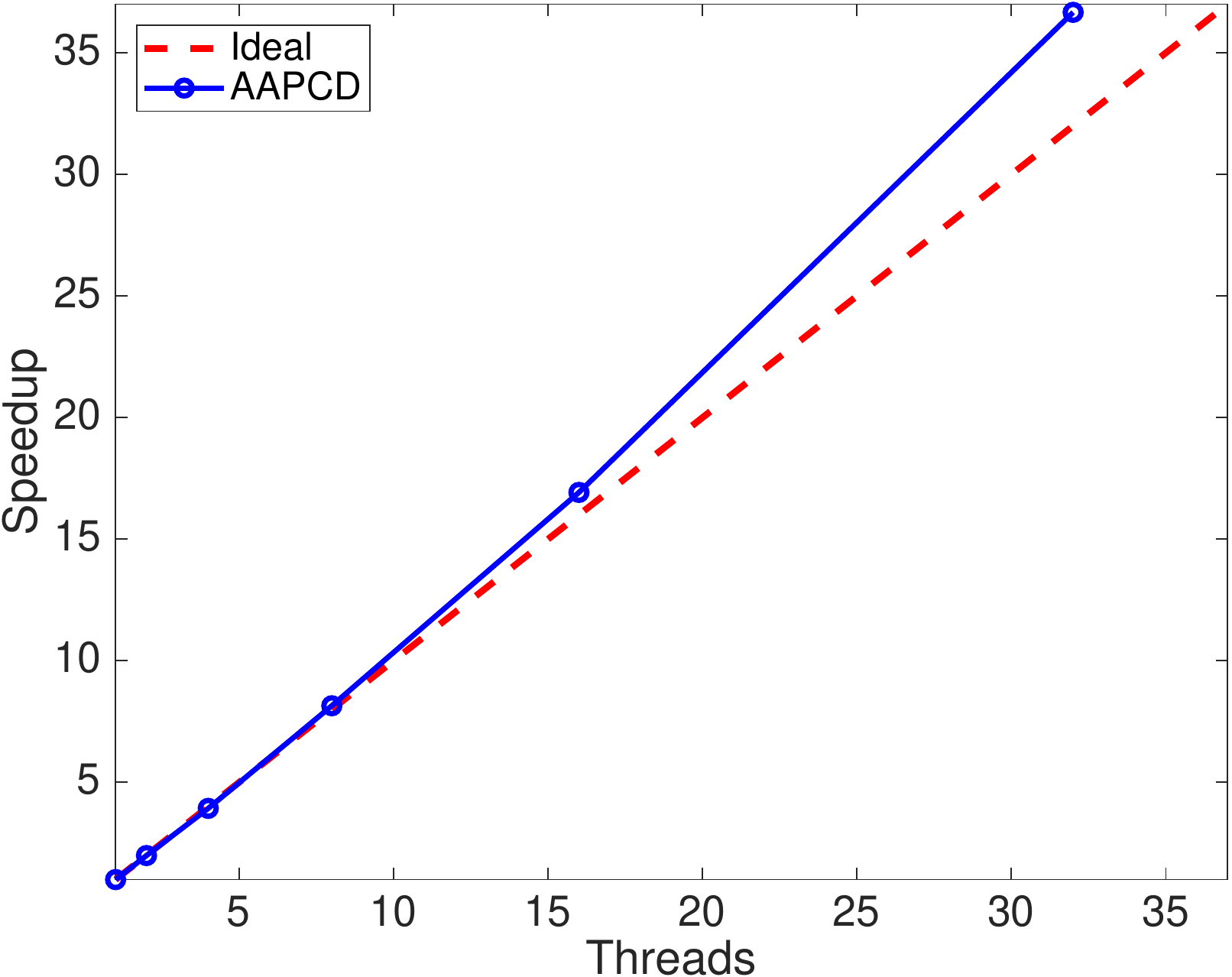}%
\caption{Speedup results of AAPCD on {\it rcv1} dataset.}
\label{fig:FSVRG_speedup}
\end{figure}

\begin{figure}[htbp]
\subfloat[a]{
\centering
\includegraphics[width=0.44\linewidth]{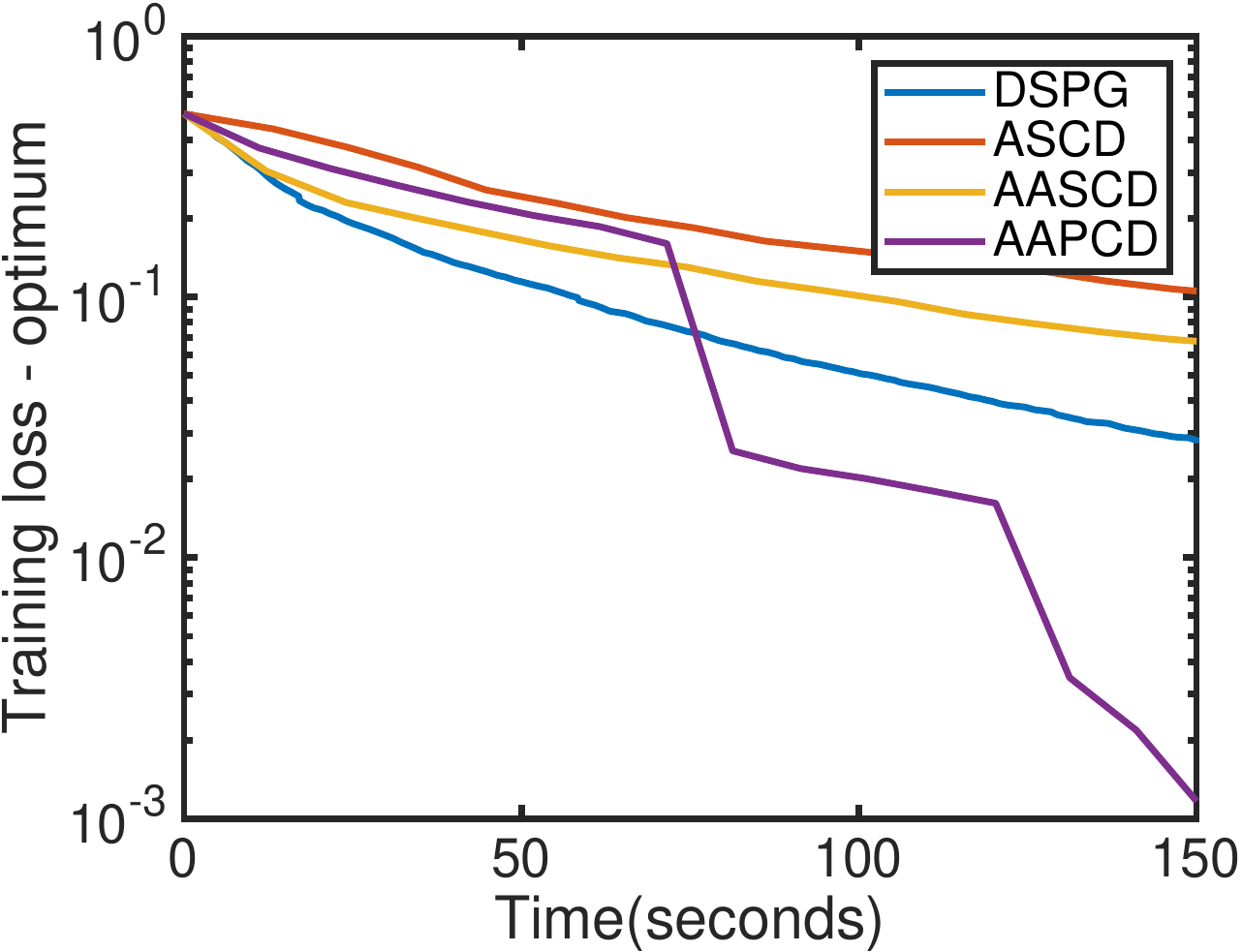}}%
\hfill
\subfloat[b]{
\centering
\includegraphics[width=0.44\linewidth]{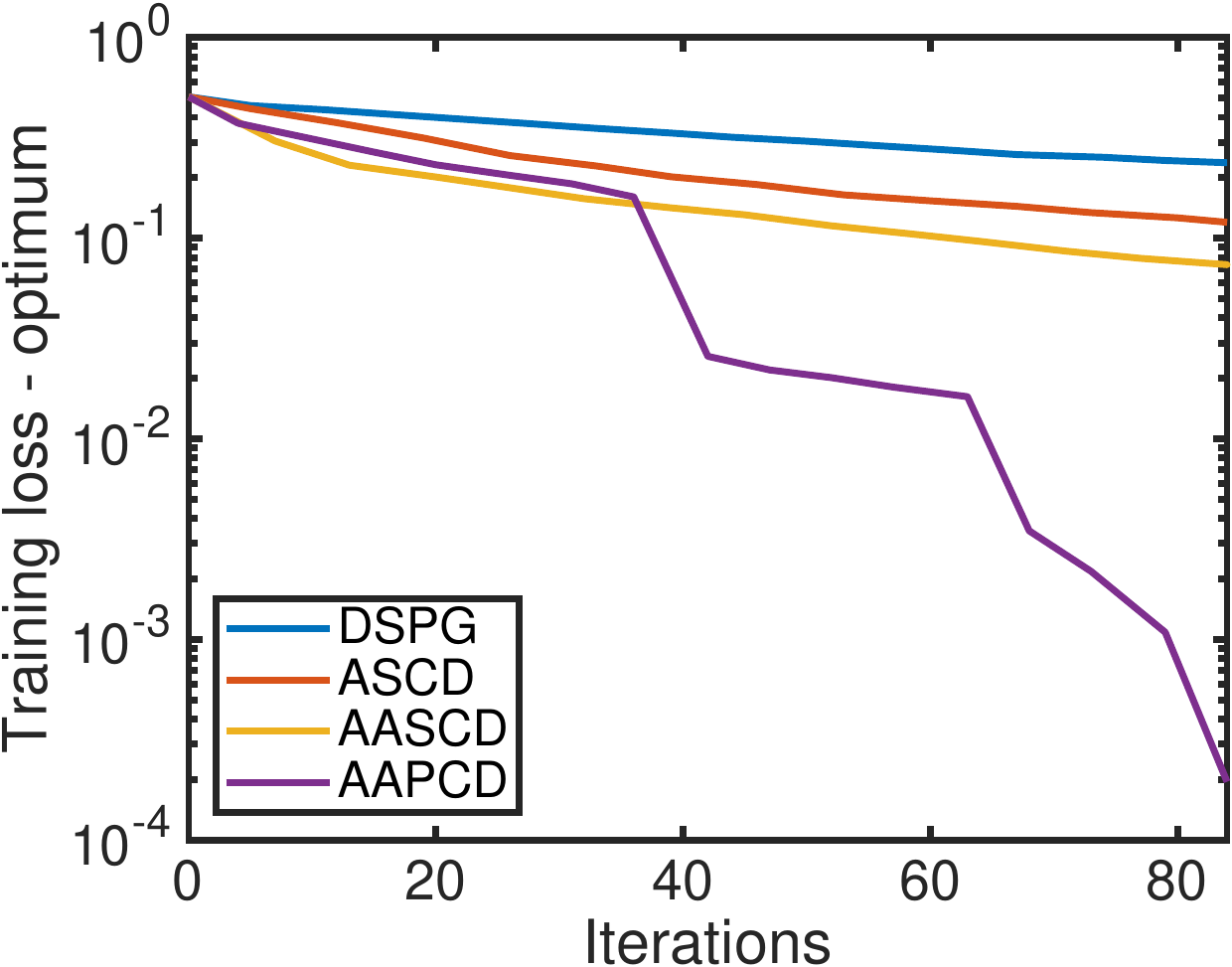}}%
\caption{Figure(a) is convergence of objective value vs. time; Figure(b) is comparison of the objective function vs. iteration for different algorithms.}
\label{fig:FSVRG_compare}
\end{figure}

In this section we test the efficiency of the asynchronous stochastic proximal coordinate descent algorithm with momentum acceleration. We performed binary classifications on the benchmark dataset {\it{rcv1}}. Following the practices in \cite{gong2013general}, we consider the logistic loss function with nonconvex regularization, 
\[
g(x) = \lambda \sum_{j=1}^m \min{(|{x_j}|,\theta)}, 
\] 
with $\lambda = 0.0001$, $\theta= 0.1\lambda$ and the zero vector as starting point. Figure \ref{fig:FSVRG_speedup} demonstrates the speedups of our algorithm. AAPCD has significant linear speedup on a parallel platform with shared memory compared to its sequential counterpart. 
We conduct experiments for comparing AAPCD with other asynchronous algorithms: ASCD \cite{liu2015asynchronous}, an synchronous version of doubly stochastic proximal algorithm (DSPG) \cite{zhao2014accelerated}, AASCD \cite{fang2018accelerating}. ASCD and DSPG did not utilize the momentum acceleration techniques. AASCD is an asynchronous accelerated variant of ASCD but only for convex and strongly convex functions. For all experiments we set the number of local workers to $32$. We set $\lambda = 0.0001$, $\theta=0.1\lambda$. For AAPCD, we set $\eta = 0.08$, $\beta=-0.08$ for negative momentum, $\beta=0.8$ for positive momentum and threshold $T_1=0.9\tau$. All blocks are of size $1000$. We set the stepsize for ASCD with $\eta=0.06$. In AASCD we set $\eta=0.09$, with momentum value $\theta_1=0.8$. For DSPG, the stepsize is $\eta=0.03$ and mini-batch size is $200$. All algorithms are terminated when the number of iterations exceeds $100$. Note that we use the best tuned parameters for each method which is obtained over a refined grid to attain the best performance. Figure \ref{fig:FSVRG_compare} shows the convergence of the objective function with respect to CPU time and the number of iterations. 
Towards the
end AAPCD decreases rapidly and needs much fewer iterations and less computing time than ASCD and AASCD to reach the same objective function values. This means that our AAPCD algorithm is very efficient and attains the best performance. Moreover AAPCD obtains a much smaller objective value by order of magnitudes compared with other algorithms. For saving space, we leave another experiment for Sigmoid loss in the supplementary materials.

Figure \ref{fig:AAPCD_neg} shows AAPCD by only applying nonnegative momentum values which is slower than AAPCD, showing that linear extrapolation using negative momentum $\beta$ for large delays is significantly useful.
\begin{figure}[htbp]
\centering
\includegraphics[width=0.6\linewidth]{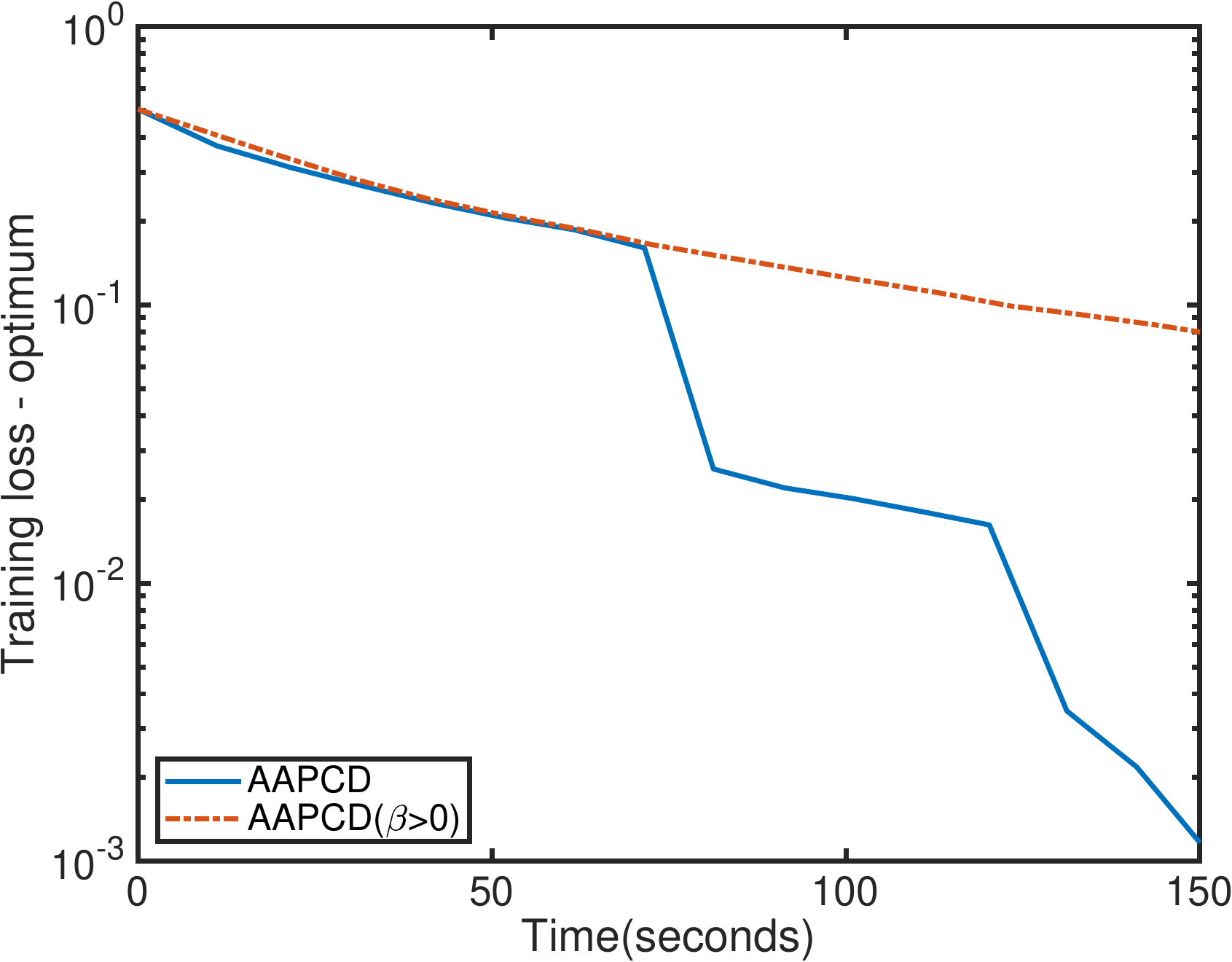}%
\caption{AAPCD versus AAPCD with momentum values $\beta >0$.}
\label{fig:AAPCD_neg}
\end{figure}

In summary our experimental results validate that AAPCD can indeed accelerate the convergence in practice.
\section{Conclusion}
In this paper, we have studied the stochastic and deterministic asynchronous parallelization of coordinate descent algorithm with momentum acceleration for efficiently solving nonconvex nonsmooth problems.
We have shown that every limit point is a critical point and proved the convergence rates for deterministic AAPCD with bounded delay. 
We verified the advantages of our method through numerical experiments.

Overall speaking, these asynchronous proximal
algorithms can be highly efficient when being used to solve
large scale nonconvex nonsmooth problems. As for future work, an extension of this study might develop the analysis in this paper to inexact proximal methods. We also plan to investigate the asynchronous parallelization of more algorithms for nonconvex nonsmooth programming for solving more complicated models.
\bibliography{main}
\bibliographystyle{aaai}
\onecolumn 

\include{sub1}

\end{document}

%% file: sub1.tex
\setcounter{equation}{0}
\setcounter{figure}{0}
\setcounter{table}{0}
\setcounter{page}{1}
\begin{center}
\textbf{\large Supplemental Materials}
\end{center}
\input{Theorems}

%% file: Theorems.tex
\noindent\subsection{Proof of Lemma \ref{P2-lemma1}}
\begin{proof}
Since $x^{k+1}_{j_k} = \Po_{j_k,\eta g_{j_k}}(y^k-\eta\nabla_{j_k} f(\hat{y}^{k}))$, we have
\begin{equation}\label{p2-lemma-eq4}
\begin{split}
\Iprod{{x}^{k+1}_{j_k}-y^{k}_{j_k}}{\nabla_{j_k} f(\hat{y}^{k})}+&\frac{1}{2\eta}\norm{{x}^{k+1}_{j_k}-y^{k}_{j_k}}^2+g_{j_k}({x}^{k+1}_{j_k})\\
&\leq g_{j_k}(y^{k}_{j_k}).
\end{split}
\end{equation}
As $f$ is $L$-Lipschitz smooth,
\[
f(x^{k+1})\leq f(y^{k}) + \Iprod{x_{j_k}^{k+1}-y_{j_k}^{k}}{\nabla_{j_k} f(y^{k})}+\frac{L}{2}\norm{x_{j_k}^{k+1}-y_{j_k}^{k}}^2.
\]
Combining with \eqref{p2-lemma-eq4}, we obtain 
\begin{equation}
\begin{split}
f(&x^{k+1})+\sum_{j=1}^m g_j(x^{k+1}_j)\leq f(y^{k})+\sum_{j=1}^m g_j(y^{k}_j)\\
&+\frac{L}{2}\norm{x^{k+1}-y^{k}}^2+\Iprod{{x}^{k+1}_{j_k}-y^{k}_{j_k}}{\nabla_{j_k} f(y^{k})-\nabla_{j_k} f(\hat{y}^{k})}-\frac{1}{2\eta}\norm{x^{k+1}-y^{k}}^2
\end{split}
\end{equation}
where we used $x^{k+1}_j=y^{k}_j$ for $j\neq j_k$. This is equivalent to,
\begin{equation}
\begin{split}
F(x^{k+1})& \leq F(y^{k})+\frac{L}{2}\norm{x^{k+1}-y^{k}}^2\\
&+\Iprod{{x}_{j_k}^{k+1}-y_{j_k}^{k}}{\nabla_{j_k} f(y^{k})-\nabla_{j_k} f(\hat{y}^{k})}-\frac{1}{2\eta}\norm{x^{k+1}-y^{k}}^2.
\end{split}
\end{equation}
For the cross term we have
\begin{equation}
\begin{split}
&\Iprod{{x}^{k+1}_{j_k}-y^{k}_{i_k}}{\nabla_{j_k} f(y^{k})-\nabla_{j_k} f(\hat{y}^{k})}\\
&\myeq{a)}{\leq} L \norm{x_{j_k}^{k+1}-y_{j_k}^k}\norm{y^k-\hat{y}^{k}}\\
&\myeq{b)}{\leq} \frac{L^2}{2C}\norm{y^k-\hat{y}^{k}}^2+\frac{C}{2}\norm{x^{k+1}-y^k}^2\\
\end{split}
\end{equation}
where $a)$ is by the Lipschitz of $\nabla f$, $b)$ by the Cauchy-Schwarz inequality. By taking expectation over $d_k$, the following sequence of inequalities is true for any $C > 0$:
\begin{equation}\label{stc-bd-deay-lemm1-eq1}
\begin{split}
&\E_{d_k}\Iprod{{x}^{k+1}_{j_k}-y^{k}_{j_k}}{\nabla_{j_k} f(y^{k})-\nabla_{i_k} f(\hat{y}^{k})~\left|~\mathcal{F}_k\right.}\\
&\leq \frac{L^2}{2C}\E_{d_k}\left[\norm{y^k-\hat{y}^{k}}^2~\left|~\mathcal{F}_k\right.\right]+\frac{C}{2}\E_{d_k}\norm{x^{k+1}-y^k}^2\\
&\myeq{a)}{\leq} \frac{L^2\tau}{2C} \sum_{h=k-\tau+1}^k\norm{y^h-y^{h-1}}^2+\frac{C}{2}\norm{x^{k+1}-y^k}^2\\
&= \left(\frac{L^2\tau}{2C} \sum_{h=k-\tau+1}^k(h-k+\tau)\norm{y^h-y^{h-1}}^2-\frac{L^2\tau}{2C} \sum_{h=k+1-\tau+1}^{k+1}(h-(k+1)+\tau)\norm{y^h-y^{h-1}}^2\right)\\
&~~~~~~~+\frac{L^2\tau^2}{2C}\norm{y^{k+1}-y^{k}}^2+\frac{C}{2}\norm{x^{k+1}-y^k}^2
\end{split}
\end{equation}
where $a)$ is due the triangle inequality and $d_k\leq \tau$. 
The linear extrapolation step for the
momentum acceleration in Algorithm \ref{AAPCD-Algo} yields  
\begin{equation}\label{Ex-Acceleration-InE}
\begin{split}
\E\norm{y^{h+1}-y^{h}}&= (1+\beta_h) \E\norm{x^{h+1}-y^{h}}\qquad \text{for } h\in {\Gamma^c}_{k}^{r}, \beta_h>0\\
\E\norm{y^{h+1}-y^{h}}&= \E\norm{x^{h+1}-y^{h}}\qquad \text{for } h\in {\Gamma^0}_{k}^{r}\\
\E\norm{y^{h+1}-y^{h}}&= (1+\beta_h) \E\norm{x^{h+1}-y^{h}}\qquad \text{for } h\in {\Gamma}_{k}^{r}, \beta_h<0.
\end{split} 
\end{equation}
Thus, by taking total expectation on both sides of \eqref{stc-bd-deay-lemm1-eq1} we have
\begin{equation}
\begin{split}
\E[\xi_{k+1}]&\leq\E[\xi_k]-\E\Iprod{{x}_{j_k}^{k+1}-y_{j_k}^{k}}{\nabla_{j_k} f(y^{k})-\nabla_{j_k} f(\hat{y}^{k})}\\
&~~~~~~~+\frac{L^2\tau^2}{2C}\E\norm{y^{k+1}-y^{k}}^2+\frac{C}{2}\E\norm{x^{k+1}-y^k}^2\\
&\leq\E[\xi_k]-\E\Iprod{{x}_{j_k}^{k+1}-y_{j_k}^{k}}{\nabla_{j_k} f(y^{k})-\nabla_{j_k} f(\hat{y}^{k})}\\
&~~~~~~~+\left(\frac{L^2\tau^2(1+\beta_k)^2}{2C}+\frac{C}{2}\right)\E\norm{x^{k+1}-y^k}^2.\\
\end{split}
\end{equation}
Therefore, we have
\begin{equation}
\begin{split}
\E[F(x^{k+1})+\xi_{k+1}]+&\frac{1}{2\eta}\E\norm{x^{k+1}-y^{k}}^2\\
\leq \E[F(y^{k})+\xi_{k}]+&\frac{L}{2}\E\norm{x^{k+1}-y^{k}}^2+\left(\frac{L^2\tau^2(1+\beta_k)^2}{2C}+\frac{C}{2}\right)\E\norm{x^{k+1}-y^k}^2.\\
\end{split}
\end{equation}
Hence, we can derive
\begin{equation}
\begin{split}
\E[F(x^{k+1})+\xi_{k+1}&]\\
\leq \E[F(y^{k})+&\xi_{k}]-\left(\frac{1}{2\eta}-\frac{L}{2}-\frac{L^2\tau^2(1+\beta_k)^2}{2C}-\frac{C}{2}\right)\E\norm{x^{k+1}-y^k}^2.\\
\end{split}
\end{equation}
By choosing $C={L\tau(1+\beta)}$, we obtain
\begin{equation}
\begin{split}
\E[F(x^{k+1})+\xi_{k+1}&]\\
\leq \E[F(y^{k})+&\xi_{k}]-\left(\frac{1}{2\eta}-\frac{L}{2}-{L\tau(1+\beta_k)}\right)\E\norm{x^{k+1}-y^k}^2\\
\end{split}
\end{equation}
and the result follows from the definition of $G(x^{k+1})$.
\end{proof}
\noindent\subsection{Proof of Theorem \ref{P2-theo1}}
\begin{proof}
Applying Lemma \ref{P2-lemma1}, we obtain that 
\begin{equation}\label{P2-Theo1-eq0}
\begin{split}
\E[G(x^{k+1})]\leq \E[G(y^{k})]-\left(\frac{1}{2\eta}-\frac{L}{2}-{L\tau(1+\beta_k)}\right)\E\norm{x^{k+1}-y^k}^2.\\
\end{split}
\end{equation}
Since $\eta < \frac{1}{L+2LT_1(1+\beta)}$ and $-1 < \beta_k < \frac{1}{L\tau}(\frac{1}{2\eta}-\frac{L}{2})-1$, it follows that $\E[G (x^{k+1})]\leq \E[G(y^k)]$.
 Moreover, the update rule of AAPCD guarantees that $F (y
^{k+1})\leq F(x^{k+1})$. In
summary, for all $k$ the following inequality holds:
\begin{equation}\label{P2-Theo1-eq1}
\E[G(y^{k+1})]\leq \E[G(x^{k+1})]\leq \E[G(y^{k})]\leq \E[G(x^{k})].
\end{equation}

Hence, from \eqref{P2-Theo1-eq0}, we obtain
\begin{equation}\label{P2-Theo1-eq01}
\begin{split}
\E[G(y^{k+1})]\leq \E[G(y^{k})]-\left(\frac{1}{2\eta}-\frac{L}{2}-{L\tau(1+\beta_k)}\right)\E\norm{x^{k+1}-y^k}^2.\\
\end{split}
\end{equation}
 From \eqref{P2-Theo1-eq01} it is seen that $\E\norm{x^{k+1}-y^k}^2$ is summable (telescoping sum). Thus, $\E\norm{x^{k+1}-y^k}\to 0$ and consequently using \eqref{Ex-Acceleration-InE} we have $\E\norm{y^{k+1}-y^k}\to 0$, which means $\E[\xi_k]\to 0$.
Combing further \eqref{P2-Theo1-eq1} with the fact that $F(x^k), F(y^k) \geq \inf F > -\infty$ for all $k$ and $\E\norm{y^{k+1}-y^k}\to 0$, we conclude that $\{G(x^k )\}, \{G(y^k)\}$ converge to the same limit $F^*$, i.e.,
\begin{equation}\label{P2-Theo1-eq2}
\lim_{k\to\infty}\E[F(x^k)] = \lim_{k\to\infty}\E[F(y^k)] =\lim_{k\to\infty}\E[G(x^k)] = \lim_{k\to\infty}\E[G(y^k)]= F^*.
\end{equation}
On the other hand, by induction we conclude from equation \eqref{P2-Theo1-eq1} that for all $k$
\[
\E[F(y^k)]\leq\E[G(y^k)]\leq F(x^0),\qquad \E[F(x^k)]\leq\E[G(x^k)]\leq F(x^0).
\]
Combining with Assumption \ref{P2-assu1} that $F$ has bounded sublevel set, we conclude that $\{x^k\}$ and $\{y^k\}$ are almost surely bounded and thus
have bounded limit points.

\noindent Since $\E\norm{x^{k+1} - y^k} \to 0$, and hence ${x^{k+1}}$ and ${y^k}$ share the same set of limit points denoted by $\Omega$. We let $l(k, j) \in \mathbb{N}$ be the last time coordinate $j$ was updated:
\[
l(k,j) = \max(\{q~|~j_q=j,q<k\}\cup \{0\}).
\]
On the other hand, by optimality condition of the proximal gradient step of AAPCD, we obtain that
\begin{equation}\label{P2-Theo1-eq41}
\begin{split}
-\nabla_{j_k} f(\hat{y}^{k})-\frac{1}{\eta}(x_{j_k}^{k+1}-y_{j_k}^k)&\in \partial g_{j_k}(x^{k+1}_{j_k})\\
\leftrightarrow\underbrace{\nabla_{j_k} f(x^{k+1})-\nabla_{j_k} f(\hat{y}^{k})-\frac{1}{\eta}(x^{k+1}_{j_k}-y^k_{j_k})}_{q^{k+1}_{j_k}}&\in \partial_{x_{j_k}} F(x^{k+1})
\end{split}
\end{equation}
and for $j\neq j_k$,
\begin{equation}\label{P2-Theo1-eq43}
\begin{split}
-\nabla_{j} f(\hat{y}^{l(k,j)})-\frac{1}{\eta}(x_{j}^{l(k,j)+1}-y_{j}^{l(k,j)})&\in \partial g_{j}(x^{k+1}_{j})\\
\leftrightarrow\underbrace{\nabla_{j} f(x^{k+1})-\nabla_{j} f(\hat{y}^{l(k,j)})-\frac{1}{\eta}(x^{l(k,j)+1}_{j}-y^{l(k,j)}_{j})}_{q^{k+1}_{j}}&\in \partial_{x_{j}} F(x^{k+1}).
\end{split}
\end{equation}
We have
\begin{equation}\label{P2-Theo1-eq5}
\begin{split}
\norm{q^{k+1}_j} &\myeq{a)}{\leq} \norm{\nabla_j f(x^{k+1})-\nabla_j f(\hat{y}^{{l(k+1,j)}})-\frac{1}{\eta}(x^{k+1}_j-y^{l(k+1,j)}_j)}\\
&\myeq{b)}{\leq} \frac{1}{\eta} \norm{x^{k+1}_j-y^{{l(k+1,j)}}_j}+L\norm{x^{k+1}_j-\hat{y}^{{l(k+1,j)}}_j}\\
&\myeq{c)}{\leq}(\frac{1}{\eta}+L)\norm{x^{{k+1}}_j-y^{{l(k+1,j)}}_j}+L\norm{y^{l(k+1,j)}_j-\hat{y}^{{l(k+1,j)}}_j}\\
&\myeq{d)}{\leq}(\frac{1}{\eta}+L)\norm{x^{{k+1}}_j-y^{{l(k+1,j)}}_j}+L\sum_{h=k+1-\tau+1-K}^{k+1}\norm{y^h_j-y^{h-1}_j}
\end{split}
\end{equation}
where $a)$ is from \eqref{P2-Theo1-eq41} and \eqref{P2-Theo1-eq43}, $b)$ by by Lipschitz of $\nabla f$, $c)$ by applying the triangle inequality and $d)$ from Assumption \ref{unbdd-determinstic-rule} and the assumption of bounded staleness. The right hand terms converge to 0 and hence $\E\norm{q^{k+1}_j}\to 0$. Since $\E\norm{x^{k+1}-y^k}^2$ is summable, it implies that $\E\norm{x^{k+1}-y^k}^2 = o(\frac{1}{k+1})$. Therefore we have $\E[\text{dist}_{\partial F(y^{k+1})}(0)]\leq \E\norm{q^{k+1}_j}=o(\frac{1}{\sqrt{k+1}})$. Consider any limit point $x\in\Omega$, and subsequences say $x^{k_t} \to x$, $y^{k_t} \to x$. By the definition
of the proximal map, the proximal gradient step of AAPCD implies that
\begin{equation}
\begin{split}
&\Iprod{\nabla_j f(\hat{y}^{l(k_t+1,j)})}{x^{k_t+1}_j-y^{l(k_t+1,j)}_j}+\frac{1}{2\eta}\norm{x^{k_t+1}_j-y^{{l({k_t+1},j)}}_j}^2+g_j(x^{k_t+1}_j)\\
&\leq\Iprod{\nabla_j f(\hat{y}^{l({k_t+1},j)})}{x_j-y^{l({k_t+1},j)}_j} + \frac{1}{2\eta}\norm{x_j-y^{l({k_t+1},j)}_j}^2+g_j(x_j).
\end{split}
\end{equation}
 Taking $\lim\sup$ on both sides and note that $x^{k+1} - y^k \to 0$, $y^{k_t} \to x$, we obtain that $\lim\sup_{t\to\infty} g_j(x^{k_t+1}_j) \leq g_j(x_j)$. Since
$g$ is lower semicontinuous and $x^{k+1}\to x$, it follows that $\lim\sup_{t\to\infty} g_j(x^{k_t+1}_j)\geq g_j(x_j)$. 
Combining both inequalities, we
conclude that $\lim_{t\to\infty} g_j(x^{k_t+1}_j) = g_j(x_j)$. Note that the continuity of $f$ yields $\lim_{t\to\infty} f(x^{k_t}) = f(x)$, we then conclude that
$\lim_{t\to\infty} F(x^{k_t}) = F(x)$, and $\lim_{t\to\infty} G(x^{k_t}) = F(x)$, since $\xi_k\to 0$. By \eqref{P2-Theo1-eq2} we have $\lim_{t\to\infty} F(x^{k_t}) = F^*$, hence
\begin{equation}\label{P2-limit-points}
 F(x)=F^*,\qquad \text{for all } x\in\Omega.
\end{equation}
Thus by \eqref{P2-limit-points}, $F$ remains constant on the compact set $\Omega$ (the set $\Omega$ is closed and bounded in $\R^m$). To this end, we have shown $x^{k_t}\to x$, $G (x^{k_t})\to F^*$. Further, we proved $q_{k_t} \in \partial F(x^{k_t})$ converges zero. We conclude that $0\in\partial F(x)$  for all $x\in\Omega$.
\end{proof}
\noindent\subsection{Proof of Theorem \ref{P2-theo2}}
\begin{proof}
Throughout the proof we assume that $r_k\neq 0$ for all $k$ because otherwise the conclusions
hold trivially. From \eqref{P2-Theo1-eq01} we have
\begin{equation}
\begin{split}
\E[G(y^{k+1})]\leq \E[G(y^{k})]-\left(\frac{1}{2\eta}-\frac{L}{2}-{L\tau(1+\beta)}\right)\E\norm{x^{k+1}-y^k}^2.\\
\end{split}
\end{equation}
By summing this inequality over $h=h-\tau-K,\ldots,k$ iterations we obtain
\begin{equation}\label{P2-Theo2-eq4}
\begin{split}
\left(\frac{1}{2\eta}-\frac{L}{2}-{L\tau(1+\beta)}\right)\sum_{h=k-K-\tau}^k\E\norm{x^{h+1}-y^{h}}^2 \leq \E[G(y^{k-K-\tau})]-\E[G(y^{k+1})].
\end{split}
\end{equation}
Moreover, equations \eqref{P2-Theo1-eq41} and \eqref{P2-Theo1-eq5} imply that
\begin{equation}\label{P2-Theo2-eq5}
\begin{split}
\left(\text{dist}_{\partial F(x^{k+1})}(0)\right)^2 &\leq \sum_{j=1}^m \left(\norm{\nabla_j f(x^{k+1})-\nabla_j f(\hat{y}^{{l(k+1,j)}})-\frac{1}{\eta}(x^{k+1}_j-y^{l(k+1,j)}_j)}\right)^2\\
&\myeq{a)}{\leq} \sum_{j=1}^m \left(\frac{1}{\eta} \norm{x^{k+1}_j-y^{{l(k+1,j)}}_j}+L\norm{x^{k+1}_j-\hat{y}^{{l(k+1,j)}}_j}\right)^2\\
&\myeq{b)}{\leq} 2(\frac{1}{\eta}+L)^2\sum_{j=1}^m\norm{x^{{k+1}}_j-y^{{l(k+1,j)}}_j}^2+2L^2\sum_{j=1}^m\norm{y^{l(k+1,j)}_j-\hat{y}^{{l(k+1,j)}}_j}^2\\
&\myeq{c)}{\leq} 2(\frac{1}{\eta}+L)^2(K+1)\sum_{h=k+1-{K}}^{k+1}\norm{x^h-y^{h-1}}^2\\
&~~~+2L^2T_1(1+\beta)\sum_{\substack{ k\in {\Gamma^c}_{k-\tau+1-{K}}^{k+1}}}\norm{x^h-y^{h-1}}^2\\
&~~~+2L^2T(1+\beta'')\sum_{\substack{ k\in {\Gamma^0}_{k-\tau+1-{K}}^{k+1}\cup{\Gamma}_{k-\tau+1-{K}}^{k+1}}}\norm{x^h-y^{h-1}}^2
\end{split}
\end{equation}
where $a)$ follows from Lipschitz of $\nabla f$, $b)$ is by the triangle inequality, $c)$ from Assumption \ref{unbdd-determinstic-rule}, the assumption of bounded delays and by \eqref{Ex-Acceleration-InE}.
We have shown in Theorem \ref{P2-theo1} that $F(x^k)\to F^*$, and it is also clear that $\text{dist}_{\Omega}(x^k)\to 0$. Thus, for any $\epsilon,\gamma > 0$ there is $k_0$ such that for all
$k \geq k_0$, we have
\begin{equation}
x^k\in\{x|\,\,\text{dist}_{\Omega}(x)\leq \epsilon, F^* < F(x) < F^*+\gamma\}.
\end{equation}
Since $\Omega$ is compact and $F$ is constant on it, the Uniformized KL property implies that for all $k \geq k_0$ 
\begin{equation}\label{P2-Theo2-eq3}
\phi'(F(x^{k+1})-F^*)\text{dist}_{\partial F(x^{k+1})}(0)\geq 1.
\end{equation}
Recall that $r_{k+1} := F(y^{k+1}) - F^*$. Then equation \eqref{P2-Theo2-eq3} is equivalent to
\begin{equation}
\begin{split}
1&\leq \left(\phi'(F(x^{k+1}) - F^*)\text{dist}_{\partial F(x^{k+1})}(0)\right)^2\\
&\leq \left(\phi'(r_{k+1})\text{dist}_{\partial F(x^{k+1})}(0)\right)^2\\
\end{split}
\end{equation}
where in the last inequality we used $r_{k+1}\leq F(x^{k+1})-F^*$ and $\phi'(t) = et^{\theta-1}$ is nonincreasing. By taking expectation on both sides of this equation and using \eqref{P2-Theo2-eq5}, we obtain

\begin{equation}
\begin{split}
\frac{1}{e^2}\E [r_{k+1}^{2-2\theta}]&{\leq} 2(\frac{1}{\eta}+L)^2(K+1)\sum_{h=k+1-{K}}^{k+1}\E\norm{x^h-y^{h-1}}^2\\
&~~~+2L^2T_1(1+\beta)\sum_{\substack{ k\in {\Gamma^c}_{k-\tau+1-K}^{k+1}}}\E\norm{x^h-y^{h-1}}^2\\
&~~~+2L^2T(1+\beta'')\sum_{\substack{ k\in {\Gamma^0}_{k-\tau+1-{K}}^{k+1}\cup{\Gamma}_{k-\tau+1-{K}}^{k+1}}}\E\norm{x^h-y^{h-1}}^2\\
&\leq \frac{2(\frac{1}{\eta}+L)^2(K+1)+2L^2T_1(1+\beta)+2L^2T(1+\beta'')}{\left(\frac{1}{2\eta}-\frac{L}{2}-{{L\tau(1+\beta)}}\right)}\\
&~~~\times\E[G(y^{k-K-\tau})-G(y^{k+1})]\\
&= b_1\E[G(y^{k-K-\tau})-G(y^{k+1})].
\end{split}
\end{equation}
The second inequality is by \eqref{P2-Theo2-eq4}, and the equality is from the definition of $b_1$. Thus by using $F(x^{k-K-\tau})\leq F(x^0)$, we have
 
\begin{equation}\label{p2-rate-relation-LK}
\E [r_{k+1}^{2-2\theta}]\leq b_1e^2 \E[G(y^{k-K-\tau})-G(y^{k+1})]\leq b_1e^2 \E[r_{0}-r_{k+1}].
\end{equation}

\noindent {\bf Part} $1:$ Suppose that $\theta = 1$, then for all $k$, we have $\E[r_0-r_{k+1}]\geq \frac{1}{b_1e^2}>0$, which cannot hold because $r_0 < \frac{1}{b_1e^2}$. Thus, $\{G(y^k)\}$ must converge in finitely many steps, which is by Theorem \ref{P2-theo1} is the stationary point of $F$.

\noindent {\bf Part} $2:$ Suppose that $\theta = \frac{1}{2}$. We have  
\begin{equation}
\E[r_{k+1}]\leq \frac{b_1e^2}{1+b_1e^2} r_{0}
\end{equation}
which yields the result.
\end{proof}
\noindent\subsection{Proof of Lemma \ref{stc-unbdd-P2-lemma1}}
\begin{proof}
Since $x^{k+1}_{j_k} = \Po_{j_k,\eta g_{j_k}}(y^k-\eta\nabla_{j_k} f(\hat{y}^{k}))$, we have
\begin{equation}\label{stc-p2-lemma-eq4}
\begin{split}
\Iprod{{x}^{k+1}_{j_k}-y^{k}_{j_k}}{\nabla_{j_k} f(\hat{y}^{k})}+&\frac{1}{2\eta}\norm{{x}^{k+1}_{j_k}-y^{k}_{j_k}}^2+g_{j_k}({x}^{k+1}_{j_k})\\
&\leq g_{j_k}(y^{k}_{j_k}).
\end{split}
\end{equation}
As $f$ is $L$-Lipschitz smooth,
\[
f(x^{k+1})\leq f(y^{k}) + \Iprod{x_{j_k}^{k+1}-y_{j_k}^{k}}{\nabla_{j_k} f(y^{k})}+\frac{L}{2}\norm{x_{j_k}^{k+1}-y_{j_k}^{k}}^2.
\]
Combining with \eqref{stc-p2-lemma-eq4}, we obtain 
\begin{equation}
\begin{split}
f(x^{k+1})&+\sum_{j=1}^m g_j(x^{k+1}_j)\leq f(y^{k})+\sum_{j=1}^m g_j(y^{k}_j)\\
&+\frac{L}{2}\norm{x^{k+1}-y^{k}}^2+\Iprod{{x}^{k+1}_{j_k}-y^{k}_{j_k}}{\nabla_{j_k} f(y^{k})-\nabla_{j_k} f(\hat{y}^{k})}-\frac{1}{2\eta}\norm{x^{k+1}_{j_k}-y^{k}_{j_k}}^2
\end{split}
\end{equation}
where we used $x^{k+1}_j=y^{k}_j$ for $j\neq j_k$. 
Therefore,
\begin{equation}\label{stc-combin-ineq}
\begin{split}
F(x^{k+1})& \leq F(y^{k})+\frac{L}{2}\norm{x^{k+1}-y^{k}}^2\\
&+\Iprod{{x}_{j_k}^{k+1}-y_{j_k}^{k}}{\nabla_{j_k} f(y^{k})-\nabla_{j_k} f(\hat{y}^{k})}-\frac{1}{2\eta}\norm{x^{k+1}-y^{k}}^2.
\end{split}
\end{equation}
Hence, 
\begin{equation}
\begin{split}
&\Iprod{{x}^{k+1}_{j_k}-y^{k}_{i_k}}{\nabla_{j_k} f(y^{k})-\nabla_{j_k} f(\hat{y}^{k})}\\
&\myeq{a)}{\leq} L \norm{x_{j_k}^{k+1}-y_{j_k}^k}\norm{y^k-\hat{y}^{k}}\\
&\myeq{b)}{\leq} \frac{L^2}{2C}\norm{y^k-\hat{y}^{k}}^2+\frac{C}{2}\norm{x^{k+1}-y^k}^2\\
\end{split}
\end{equation}
where $a)$ is by the Lipschitz of $\nabla f$ and $b)$ is by Cauchy-Schwarz inequality. We bound the expectation of $\norm{y^k-\hat{y}^{k}}^2$ over the delay. In particular, the following sequence of inequalities is true for any $C > 0$:
\begin{equation}\label{stc-combin-ineq-cross}
\begin{split}
&\E_{d_k}\Iprod{{x}^{k+1}_{j_k}-y^{k}_{j_k}}{\nabla_{j_k} f(y^{k})-\nabla_{i_k} f(\hat{y}^{k})~\left|~\mathcal{F}_k\right.}\\
&\leq \frac{L^2}{2C}\E_{d_k}\left[\norm{y^k-\hat{y}^{k}}^2~\left|~\mathcal{F}_k\right.\right]+\frac{C}{2}\E_{d_k}\norm{x^{k+1}-y^k}^2\\
&\leq \frac{L^2}{2C} \E_{d_k}\sum_{j=1}^mi\sum_{h=k-i+1}^k\norm{y^h_j-y^{h-1}_j}^2+\frac{C}{2}\norm{x^{k+1}-y^k}^2\\
&\leq \frac{L^2}{2C} \sum_{i=1}^{\infty}ip_{i}\sum_{h=k-i+1}^k\norm{y^h-y^{h-1}}^2+\frac{C}{2} \norm{x^{k+1}-y^k}^2\\
&\myeq{a)}{=} \left(\frac{L^2}{2C} \sum_{i=1}^{\infty}ip_{i}\sum_{h=k-i+1}^k(h-k+i)\norm{y^h-y^{h-1}}^2-\frac{L^2}{2C} \sum_{i=1}^{\infty}ip_{i}\sum_{h=k-i+2}^{k+1}(h-(k+1)+i)\norm{y^h-y^{h-1}}^2\right)\\
&~~~~~~~+\frac{L^2c_0}{2C}\norm{y^{k+1}-y^{k}}^2+\frac{C}{2}\norm{x^{k+1}-y^k}^2.\\
&\myeq{b)}{=} \left(\frac{L^2}{2C} \sum_{h=1}^{k}\sum_{i=k-h+1}^\infty ip_{i}(h-k+i)\norm{y^h-y^{h-1}}^2-\frac{L^2}{2C} \sum_{h=1}^{k+1}\sum_{i=k-h+2}^{\infty}ip_{i}(h-(k+1)+i)\norm{y^h-y^{h-1}}^2\right)\\
&~~~~~~~+\frac{L^2c_0}{2C}\norm{y^{k+1}-y^{k}}^2+\frac{C}{2}\norm{x^{k+1}-y^k}^2\\
&= \left(\frac{L^2}{2C} \sum_{h=1}^{k}\sum_{t=1}^\infty t(t-h+k)p_{t-h+k}\norm{y^h-y^{h-1}}^2-\frac{L^2}{2C} \sum_{h=1}^{k+1}\sum_{t=1}^{\infty}t(t-h+k+1)p_{t-h+k+1}\norm{y^h-y^{h-1}}^2\right)\\
&~~~~~~~+\frac{L^2c_0}{2C}\norm{y^{k+1}-y^{k}}^2+\frac{C}{2}\norm{x^{k+1}-y^k}^2\\
&\myeq{c)}{=}\left(\frac{L^2}{2C} \sum_{h=1}^{k} c_{k-h}\norm{y^h-y^{h-1}}^2-\frac{L^2}{2C} \sum_{h=1}^{k+1}c_{k+1-h}\norm{y^h-y^{h-1}}^2\right)\\
&~~~~~~~+\frac{L^2c_0}{2C}\norm{y^{k+1}-y^{k}}^2+\frac{C}{2}\norm{x^{k+1}-y^k}^2
\end{split}
\end{equation}
where in $a)$, we used $c_0=\sum_{t=1}^{\infty}t^2p_t$, in $b)$, we switched the order of summation in the double sum, and c) uses $c_{k-h}:=\sum_{t=1}^\infty t(t-h+k)p_{t-h+k}$. Taking total expectation on the equation above, we obtain 
\begin{equation}
\begin{split}
&\E\Iprod{{x}^{k+1}_{j_k}-y^{k}_{j_k}}{\nabla_{j_k} f(y^{k})-\nabla_{j_k} f(\hat{y}^{k})}\\
&\leq\left(\frac{L^2}{2C} \sum_{h=1}^{k} c_{k-h}\E\norm{y^h-y^{h-1}}^2-\frac{L^2}{2C} \sum_{h=1}^{k+1}c_{k+1-h}\E\norm{y^h-y^{h-1}}^2\right)\\
&~~~~~~~+\frac{L^2c_0}{2C}\E\norm{y^{k+1}-y^{k}}^2+\frac{C}{2}\E\norm{x^{k+1}-y^k}^2\\
&\myeq{d)}{\leq} \left(\frac{L^2}{2C} \sum_{h=1}^{k} c_{k-h}\E\norm{y^h-y^{h-1}}^2-\frac{L^2}{2C} \sum_{h=1}^{k+1}c_{k+1-h}\E\norm{y^h-y^{h-1}}^2\right)\\
&~~~~~~~+\left(\frac{L^2c_0(1+\beta_k)^2}{2C}+\frac{C}{2}\right)\E\norm{x^{k+1}-y^k}^2\\
&=\E[\xi_k]-\E[\xi_{k+1}]+\left(\frac{L^2c_0(1+\beta_k)^2}{2C}+\frac{C}{2}\right)\E\norm{x^{k+1}-y^k}^2
\end{split}
\end{equation}
where $d)$ follows from \eqref{Ex-Acceleration-InE}. 
Thus, we have
\begin{equation}
\begin{split}
\E[\xi_{k+1}]&\leq\E[\xi_k]-\E\Iprod{{x}_{j_k}^{k+1}-y_{j_k}^{k}}{\nabla_{j_k} f(y^{k})-\nabla_{j_k} f(\hat{y}^{k})}\\
&~~~~~~~+\left(\frac{L^2c_0(1+\beta_k)^2}{2C}+\frac{C}{2}\right)\E\norm{x^{k+1}-y^k}^2\\
\end{split}
\end{equation}
which is substituted into \eqref{stc-combin-ineq} to yield
\begin{equation}
\begin{split}
\E [F(x^{k+1})]+&\E[\xi_{k+1}]\leq \E [F(y^{k})]+\E[\xi_{k}]-\frac{1}{2\eta}\E\norm{x^{k+1}-y^{k}}^2\\
&+\frac{L}{2}\E\norm{x^{k+1}-y^{k}}^2+\left(\frac{L^2c_0(1+\beta_k)^2}{2C}+\frac{C}{2}\right)\E\norm{x^{k+1}-y^k}^2.\\
\end{split}
\end{equation}
Thus, we get
\begin{equation}
\begin{split}
\E[F(x^{k+1})+&\xi_{k+1}]\\
\leq \E[F(y^{k})+\xi_{k}]-&\left(\frac{1}{2\eta}-\frac{L}{2}-\frac{L^2c_0(1+\beta_k)^2}{2C}-\frac{C}{2}\right)\E\norm{x^{k+1}-y^k}^2.\\
\end{split}
\end{equation}
Finally, by choosing $C=L(1+\beta)\sqrt{{c_0}}$, we have
\begin{equation}
\begin{split}
\E[F(x^{k+1})+&\xi_{k+1}]\\
\leq \E[F(y^{k})+\xi_{k}]-&\left(\frac{1}{2\eta}-\frac{L}{2}-{L(1+\beta_k){\sqrt{c_0}}}\right)\E\norm{x^{k+1}-y^k}^2.\\
\end{split}
\end{equation}
In particular for all $k\in\mathbb{N}$, we have $G(y^k) = F(y^k)+\xi_k$ so \eqref{stc-unbdd-p2-lemma-eq} follows. 
\end{proof}
\noindent\subsection{Proof of Theorem \ref{stc-unbdd-P2-theo1}}
\begin{proof}
From Lemma \ref{stc-unbdd-P2-lemma1} we have
\begin{equation}\label{stc-unbdd-P2-Theo1-eq0}
\begin{split}
\E[G(x^{k+1})]\leq \E[G(y^{k})]-\left(\frac{1}{2\eta}-\frac{L}{2}-{L(1+\beta_k){\sqrt{c_0}}}\right)\E\norm{x^{k+1}-y^k}^2.\\
\end{split}
\end{equation}
Since $\eta<\frac{1}{L+2L{\sqrt{c_{T_1}}}(1+\beta)}$ and by the upper bound for $\beta_k$, it follows that $\E[G (x^{k+1})]\leq \E[G(y^k)]$. Moreover, the update rule of AAPCD guarantees that $F (y^{k+1})\leq F(x^{k+1})$ and so we have $\E[G (y^{k+1})]\leq \E[G(x^{k+1})]$. In
summary, for all $k$ the following inequality holds:
\begin{equation}\label{stc-monoto-f-g}
\E [G(y^{k+1})]\leq \E [G(x^{k+1})]\leq \E [G(y^{k})]\leq \E [G(x^{k})].
\end{equation}
Thus, from \eqref{stc-unbdd-P2-Theo1-eq0} we obtain
\begin{equation}\label{stc-unbdd-P2-Theo1-eq01}
\begin{split}
\E[G(y^{k+1})]\leq \E[G(y^{k})]-\left(\frac{1}{2\eta}-\frac{L}{2}-{L(1+\beta_k){\sqrt{c_0}}}\right)\E\norm{x^{k+1}-y^k}^2.\\
\end{split}
\end{equation}
 Hence $\E\norm{x^{k+1}-y^{k}}^2$ is summable and $\E\norm{x^{k+1}-y^{k}}^2$ converges to zero. 
Further, by the fact $\sum_{h=1}^\infty\E\norm{ y^{h}-y^{h-1}}^2< \infty$ and using $\sum_{h=1}^\infty c_{h} < \infty$, the series $S_k = \sum_{h=1}^k c_{k-h} \E\norm{y^{h}-y^{h-1}}^2$ converges to zero, i.e., $
\E[\xi_k]\to 0$. Combing further with the fact that $F(x^k), F(y^k) \geq \inf F > -\infty$ for all $k$ and \eqref{stc-monoto-f-g}, we conclude that $\{F(x^k )\}, \{F(y^k)\}$ converge to the same limit $F^*$, i.e.,
\begin{equation}\label{stc-unbdd-P2-Theo1-eq2}
\lim_{k\to\infty}\E[F(x^k)] = \lim_{k\to\infty}\E[F(y^k)] =\lim_{k\to\infty}\E[G(x^k)] = \lim_{k\to\infty}\E[G(y^k)]= F^*.
\end{equation}
Since $\E\norm{y^k - x^{k+1}} \to 0$, and hence ${x^k}$ and ${y^k}$ share the same set of limit points which is denoted by $\Omega$. Similar to  the analysis of AAPCD with bounded delay, using the optimality condition of the proximal step of AAPCD, we obtain that
\begin{equation}\label{stc-P2-Theo1-eq4}
\begin{split}
-\nabla_{j_k} f(\hat{y}^{k})-\frac{1}{\eta}(x_{j_k}^{k+1}-y_{j_k}^k)\in \partial g_{j_k}(x^{k+1}_{j_k})&\\
\Longrightarrow\underbrace{\nabla_{j_k} f(x^{k+1})-\nabla_{j_k} f(\hat{y}^{k})-\frac{1}{\eta}(x^{k+1}_{j_k}-y^k_{j_k})}_{q^{k+1}_{j_k}}&\in \partial_{x_{j_k}} F(x^{k+1})
\end{split}
\end{equation}
and for $j\neq j_k$,
\begin{equation}\label{stc-P2-Theo1-eq41}
\begin{split}
-\nabla_{j} f(\hat{y}^{l(k,j)})-\frac{1}{\eta}(x_{j}^{l(k,j)+1}-y_{j}^{l(k,j)})\in \partial g_{j}(x^{k+1}_{j})&\\
\Longrightarrow\underbrace{\nabla_{j} f(x^{k+1})-\nabla_{j} f(\hat{y}^{l(k,j)})-\frac{1}{\eta}(x^{k+1}_{j}-y^{l(k,j)}_{j})}_{q^{k+1}_{j}}&\in \partial_{x_{j}} F(x^{k+1})
\end{split}.
\end{equation}
By the Assumption \ref{unbdd-determinstic-rule}, we can derive
\begin{equation}\label{stc-P2-Theo1-eq5}
\begin{split}
\E\norm{q^{k+1}_j}^2 &\leq \E\left(\norm{\nabla_j f(x^{k+1})-\nabla_j f(\hat{y}^{{l(k+1,j)}})-\frac{1}{\eta}(x^{k+1}_j-y^{l(k+1,j)}_j)}\right)^2\\
&\myeq{a)}{\leq} \E\left(\frac{1}{\eta} \norm{x^{k+1}_j-y^{{l(k+1,j)}}_j}+L\norm{x^{k+1}-\hat{y}^{{l(k+1,j)}}}\right)^2\\
&\myeq{b)}{\leq} \E\left((\frac{1}{\eta}+L)\norm{x^{{k+1}}_j-y^{{l(k+1,j)}}_j}+L\norm{y^{l(k+1,j)}-\hat{y}^{{l(k+1,j)}}}\right)^2\\
&{\leq} 2(\frac{1}{\eta}+L)^2\E\norm{x^{{k+1}}_j-y^{{l(k+1,j)}}_j}^2+2L^2\E\norm{y^{l(k+1,j)}-\hat{y}^{{l(k+1,j)}}}^2\\
&\myeq{c)}{\leq} 2(\frac{1}{\eta}+L)^2\E\norm{x^{{k+1}}_j-y^{{l(k+1,j)}}_j}^2+2L^2\sum_{h=k-{K}}^k\E[\xi_{h}-\xi_{h+1}]\\
&~~~+{2L^2c_0}\sum_{h=k-K}^{k}\E\norm{y^{h+1}-y^{h}}^2\\
&\myeq{d)}{\leq} 2(\frac{1}{\eta}+L)^2K\sum_{h=k-{K}}^{k}\E\norm{x^{h+1}_j-y^{h}_j}^2+2L^2\E[\xi_{k-{K}}-\xi_{k+1}]\\
&~~~+{2L^2c_0}\sum_{h=k-K}^{k}\E\norm{y^{h+1}-y^{h}}^2
\end{split}
\end{equation}
where $a)$ is by Lipschitz of $\nabla f$, $b)$ by the triangle inequality, $c)$ from \eqref{stc-combin-ineq-cross} and $d)$ by using a telescoping sum. As $k\to \infty$, this right term converges to zero and therefore we have $\E\norm{q_j^{k+1}}^2\to 0$. 
On the other hand, by induction we conclude from equation \eqref{stc-unbdd-P2-Theo1-eq2} that for all $k$
\[
\E[F(y^k)]\leq \E[G(y^k)]\leq F(x^0)\qquad \E[F(x^k)]\leq \E[G(x^k)]\leq F(x^0).
\]
Combining with Assumption \ref{P2-assu1} that $F$ has bounded sublevel set, we conclude that $\{x^k\}$ and $\{y^k\}$ are almost surely bounded and thus
have bounded limit points.
We assume that $\{x^k\}$ and $\{y^k\}$ are bounded and thus have bounded limit points.
We fix any limit point $x\in\Omega$, say $x^{k_t} \to x$, $y^{k_t} \to x$. Note that the continuity of $f$ yields $\lim_{t\to\infty} f(x^{k_t}) = f(x)$. Moreover, by the definition
of the proximal map, the proximal gradient step of AAPCD implies that
\begin{equation}
\begin{split}
&\Iprod{\nabla_j f(\hat{y}^{l(k_t+1,j)})}{x^{k_t+1}_j-y^{l(k_t+1,j)}_j}+\frac{1}{2\eta}\norm{x^{k_t+1}_j-y^{l(k_t+1,j)}_j}^2+g_j(x^{k_t+1}_j)\\
&\leq\Iprod{\nabla_j f(\hat{y}^{l(k_t+1,j)})}{x_j-y^{l(k_t+1,j)}_j} + \frac{1}{2\eta}\norm{x_j-y^{l(k_t+1,j)}_j}^2+g_j(x_j).
\end{split}
\end{equation}
Hence, we have
\begin{equation}
\begin{split}
g_j(x^{k_t+1}_j)&\leq\Iprod{\nabla_j f(\hat{y}^{l(k_t+1,j)})}{x_j-x^{k_t+1}_j} + \frac{1}{2\eta}\norm{x_j-y^{l(k_t+1,j)}_j}^2+g_j(x_j)\\
&\leq\Iprod{\nabla_j f(\hat{y}^{l(k_t+1,j)})}{x_j-x^{k_t+1}_j}\\
&~~~ + \frac{1}{\eta}\norm{x_j - y_j^{k_t+1}}^2+\frac{1}{\eta}\norm{y_j^{k_t+1}-y^{l(k_t+1,j)}_j}^2+g_j(x_j)
\end{split}
\end{equation}
 where the last inequality is by the triangle inequality. Taking $\lim\sup$ on both sides and note that $x^{k_t} - y^{k_t} \to 0$, $y^{k_t} \to x$, we obtain that $\lim\sup_{t\to\infty} g_j(x^{{k_t}}_j) \leq g_j(x)$. Since
$g_j$ is lower semicontinuous and $x^{k_t}\to x$, it follows that $\lim\sup_{t\to\infty} g_j(x^{k_t}_j)\geq g_j(x)$. 
Combining both inequalities we
conclude that $\lim_{t\to\infty} g_j(x^{k_t}_j) = g_j(x)$. Hence we have that
$\lim_{t\to\infty} F(x^{k_t}) = F(x)$. Since $\lim_{t\to\infty} F(x^{k_t}) = F^*$ by equation \eqref{stc-unbdd-P2-Theo1-eq2}, we get
\begin{equation}
 F(x)=F^*,\qquad \text{for all } x\in\Omega.
\end{equation}
Hence, $F$ remains constant on the compact set $\Omega$. Since $\lim \E[\xi_k]\to 0$, we have, $\E[G (x^{k_t})]\to F(x)$ and thus
$\E[G (x^{k_t})]\to F^*$. We have shown $x^{k_t}\to x$, $F (x^{k_t})\to F(x)$ and that
$q_{k_t} \in \partial F(x^{k_t})$ converges zero. Altogether, we have $0\in \partial F(x)$ for all $x\in\Omega$.
\end{proof}
\noindent\subsection{Proof of Theorem \ref{stc-unbdd-P2-theo2}}
\begin{proof}
Throughout the proof we assume that $r_k\neq 0$ for all $k$ because otherwise the algorithm terminates and the conclusions
hold trivially. Lemma \ref{stc-unbdd-P2-lemma1} yields that
\begin{equation}\label{stc-P2-Theo2-eq1}
\begin{split}
\E[G(x^{h+1})]\leq \E[G(y^{h})]-\left(\frac{1}{2\eta}-\frac{L}{2}-{L(1+\beta_k){\sqrt{c_0}}}\right)\E\norm{x^{h+1}-y^h}^2.\\
\end{split}
\end{equation}
Therefore, from \eqref{stc-monoto-f-g} we obtain
\begin{equation}
\begin{split}
\E[G(y^{h+1})]\leq \E[G(y^{h})]-\left(\frac{1}{2\eta}-\frac{L}{2}-{L(1+\beta_k){\sqrt{c_0}}}\right)\E\norm{x^{h+1}-y^h}^2.\\
\end{split}
\end{equation}
By summing this inequality over $h=k-K,\ldots,k$ iterations we obtain
\begin{equation}\label{stc-theo2-eq2}
\begin{split}
\left(\frac{1}{2\eta}-\frac{L}{2}-{L(1+\beta){\sqrt{c_0}}}\right)\sum_{h=k-K}^k\E\norm{x^{h+1}-y^{h}}^2 \leq \E[G(y^{k-K})]-\E[G(y^{k+1})].
\end{split}
\end{equation}
We have shown in Theorem \ref{stc-unbdd-P2-theo1} that $F(x^k)\to F^*$, and it is also clear that $\text{dist}_{\Omega}(x^k)\to 0$. Thus, for any $\epsilon,\gamma > 0$ there is $k_0$ such that for all
$k \geq k_0$, we have
\begin{equation}
x^k\in\{x|\,\,\text{dist}_{\Omega}(x)\leq \epsilon, F^* < F(x) < F^*+\gamma\}
\end{equation}
Since $\Omega$ is compact and $F$ is constant on it, we can apply KL property. The Uniformized KL property implies that for all $k \geq k_0$ 
\begin{equation}\label{stc-P2-Theo2-eq3}
\phi'(F(x^{k+1})-F^*)\text{dist}_{\partial F(x^{k+1})}(0)\geq 1.
\end{equation}
Moreover, equations \eqref{stc-P2-Theo1-eq4} and \eqref{stc-P2-Theo1-eq41} imply that
\begin{equation}
\begin{split}
1&\myeq{a)}{\leq} \left(\phi'(F(x^{k+1})-F^*)\text{dist}_{\partial F(x^{k+1})}(0)\right)^2\\
&\myeq{b)}{\leq} \left(\phi'(r_{k+1})\text{dist}_{\partial F(x^{k+1})}(0)\right)^2\\
&\leq (\phi'(r_{k+1}))^2 \sum_{j=1}^m\left(\norm{\nabla_j f(x^{k+1})-\nabla_j f(\hat{y}^{{l(k+1,j)}})-\frac{1}{\eta}(x^{k+1}_j-y^{l(k+1,j)}_j)}\right)^2\\
&\leq (\phi'(r_{k+1}))^2 \sum_{j=1}^m\left(\norm{\nabla_j f(x^{k+1})-\nabla_j f(\hat{y}^{{l(k+1,j)}})-\frac{1}{\eta}(x^{k+1}_j-y^{l(k+1,j)}_j)}\right)^2\\
&\leq (\phi'(r_{k+1}))^2 \left(\frac{2}{\eta^2} \sum_{j=1}^m\norm{x^{k+1}_j-y^{{l(k+1,j)}}_j}^2+2\sum_{j=1}^m\norm{\nabla_j f(x^{k+1})-\nabla_j f(\hat{y}^{{l(k+1,j)}})}^2\right)\\
&\myeq{c)}{\leq} (\phi'(r_{k+1}))^2 \left(\frac{2}{\eta^2} \sum_{j=1}^m\norm{x^{k+1}_j-y^{{l(k+1,j)}}_j}^2+2L^2\sum_{j=1}^m\norm{x^{k+1}-\hat{y}^{{l(k+1,j)}}}^2\right)\\
&\myeq{d)}{\leq} (\phi'(r_{k+1}))^2 \left((\frac{2}{\eta^2}+4L^2)\sum_{j=1}^m\norm{x^{{k+1}}_j-y^{{l(k+1,j)}}_j}^2+4L^2\sum_{j=1}^m\norm{y^{l(k+1,j)}-\hat{y}^{{l(k+1,j)}}}^2\right)\\
\end{split}
\end{equation}
where $a)$ follows from \eqref{stc-P2-Theo2-eq3}, $b)$ by $r_{k+1} \leq F(x^{k+1})-F^*$ and the fact that $\phi'$ is nonincreasing, $c)$ from the Lipschitz of $\nabla f$ and  $d)$ from the triangle inequality. 
We have that $\phi'(t)={e}t^{\theta-1}$. Thus the above equation becomes
\begin{equation} 
\frac{1}{e^2}r_{k+1}^{2-2\theta} \leq (\frac{2}{\eta^2}+4L^2)\sum_{j=1}^m\norm{x^{{k+1}}_j-y^{{l(k+1,j)}}_j}^2+4L^2\sum_{j=1}^m\norm{y^{l(k+1,j)}-\hat{y}^{{l(k+1,j)}}}^2.
\end{equation}
By taking total expectation on both sides of this equation, and following the derivations  similar to that of \eqref{stc-P2-Theo1-eq5}, we have
\begin{equation}
\begin{split}
\frac{1}{e^2}\E [r_{k+1}^{2-2\theta}] &\leq (\frac{2}{\eta^2}+4L^2)\sum_{j=1}^m\E\norm{x^{{k+1}}_j-y^{{l(k+1,j)}}_j}^2+4L^2\sum_{j=1}^m\E\norm{y^{l(k+1,j)}-\hat{y}^{{l(k+1,j)}}}^2\\
&\leq (\frac{2}{\eta^2}+4L^2)\sum_{h=k-{K}}^{k}\E\norm{x^{{h+1}}-y^{h}}^2+4L^2\E[\xi_{k-{K}}-\xi_{k+1}]\\
&~~~+{4L^2c_0}\sum_{h=k-K}^{k}\E\norm{y^{h+1}-y^{h}}^2\\
&\leq (\frac{2}{\eta^2}+4L^2)\sum_{h=k-{K}}^{k}\E\norm{x^{{h+1}}-y^{h}}^2+4L^2\E[\xi_{k-{K}}-\xi_{k+1}]\\
&~~~+{4L^2c_0}(1+\beta)\sum_{\substack{ k\in {\Gamma^c}_{k+1-{K}}^{k+1}}}\norm{x^h-y^{h-1}}^2+{4L^2c_0}\sum_{\substack{ k\in {\Gamma^0}_{k+1-{K}}^{k+1}}}\norm{x^h-y^{h-1}}^2\\
&~~~+{4L^2c_0}(1+\beta')\sum_{\substack{ k\in {\Gamma}_{k+1-{K}}^{k+1}}}\norm{x^h-y^{h-1}}^2\\
&\leq b_1\E[G(y^{k-K})-G(y^{k+1})]\\
&\leq b_1\E[F(y^{0})-F(y^{k+1})]\\
&= b_1\E[r_0-r_{k+1}]
\end{split}
\end{equation}
where the second last inequality is by the definition of $b_1$ and \eqref{stc-theo2-eq2}. Thus from the above inequality we have 
\begin{equation}\label{stc-p2-rate-relation-LK}
\E [r_{k+1}^{2-2\theta}]\leq b_1e^2 \E[r_0-r_{k+1}].
\end{equation}
\noindent {\bf Part} $1:$ Suppose that $\theta = 1$, then for all $k$, we have $\E[r_0-r_{k+1}]\geq \frac{1}{b_1e^2}>0$, which cannot hold because $\E[r_{k+1}]\geq 0$ and $r_0< \frac{1}{b_1e^2}$. Thus, $\{F(y^k)\}$ must converge in finitely many steps, which is by Theorem \ref{stc-unbdd-P2-theo1} the stationary point of $F$.

\noindent {\bf Part} $2:$ Suppose that $\theta = \frac{1}{2}$. We have from \eqref{stc-p2-rate-relation-LK}
\begin{equation}
\E[r_{k+1}]\leq \frac{b_1e^2}{1+b_1e^2} \E[r_{0}]
\end{equation}
which yields the result.
\end{proof}
\noindent\subsection{Proof of Lemma \ref{unbdd-P2-lemma1}}
\begin{proof}
Since $x^{k+1}_{j_k} = \Po_{j_k,\eta g_{j_k}}(y^k-\eta\nabla_{j_k} f(\hat{y}^{k}))$, we have
\begin{equation}\label{unbd-p2-lemma-eq4}
\begin{split}
\Iprod{{x}^{k+1}_{j_k}-y^{k}_{j_k}}{\nabla_{j_k} f(\hat{y}^{k})}+&\frac{1}{2\eta}\norm{{x}^{k+1}_{j_k}-y^{k}_{j_k}}^2+g_{j_k}({x}^{k+1}_{j_k})\\
&\leq g_{j_k}(y^{k}_{j_k}).
\end{split}
\end{equation}
As $f$ is $L$-Lipschitz smooth,
\[
f(x^{k+1})\leq f(y^{k}) + \Iprod{x_{j_k}^{k+1}-y_{j_k}^{k}}{\nabla_{j_k} f(y^{k})}+\frac{L}{2}\norm{x_{j_k}^{k+1}-y_{j_k}^{k}}^2.
\]
Combining with \eqref{unbd-p2-lemma-eq4}, we obtain 
\begin{equation}
\begin{split}
f(&x^{k+1})+\sum_{j=1}^m g_j(x^{k+1}_j)\leq f(y^{k})+\sum_{j=1}^m g_j(y^{k}_j)\\
&+\frac{L}{2}\norm{x^{k+1}-y^{k}}^2+\Iprod{{x}^{k+1}_{j_k}-y^{k}_{j_k}}{\nabla_{j_k} f(y^{k})-\nabla_{j_k} f(\hat{y}^{k})}-\frac{1}{2\eta}\norm{x^{k+1}-y^{k}}^2
\end{split}
\end{equation}
where we used $x^{k+1}_j=y^{k}_j$ for $j\neq j_k$. This is equivalent to,
\begin{equation}
\begin{split}
F(x^{k+1})& \leq F(y^{k})+\frac{L}{2}\norm{x^{k+1}-y^{k}}^2\\
&+\Iprod{{x}_{j_k}^{k+1}-y_{j_k}^{k}}{\nabla_{j_k} f(y^{k})-\nabla_{j_k} f(\hat{y}^{k})}-\frac{1}{2\eta}\norm{x^{k+1}-y^{k}}^2.
\end{split}
\end{equation}
For the cross term we have
\begin{equation}
\begin{split}
&\Iprod{{x}^{k+1}_{j_k}-y^{k}_{i_k}}{\nabla_{j_k} f(y^{k})-\nabla_{j_k} f(\hat{y}^{k})}\\
&\myeq{a)}{\leq} L \norm{x_{j_k}^{k+1}-y_{j_k}^k}\norm{y^k-\hat{y}^{k}}\\
&\myeq{b)}{\leq} \frac{L^2}{2C}\norm{y^k-\hat{y}^{k}}^2+\frac{C}{2}\norm{x^{k+1}-y^k}^2\\
\end{split}
\end{equation}
where $a)$ is by the Lipschitz of $\nabla f$, $b)$ by the Cauchy-Schwarz inequality. By taking expectation over $d_k$, the following sequence of inequalities is true for any $C > 0$:
\begin{equation}\label{stc-unbd-deay-lemm1-eq1}
\begin{split}
&\E_{d_k}\Iprod{{x}^{k+1}_{j_k}-y^{k}_{j_k}}{\nabla_{j_k} f(y^{k})-\nabla_{i_k} f(\hat{y}^{k})~\left|~\mathcal{F}_k\right.}\\
&\leq \frac{L^2}{2C}\E_{d_k}\left[\norm{y^k-\hat{y}^{k}}^2~\left|~\mathcal{F}_k\right.\right]+\frac{C}{2}\E_{d_k}\norm{x^{k+1}-y^k}^2\\
&\myeq{a)}{\leq} \frac{L^2\tau}{2C} \sum_{h=k-\tau+1}^k\norm{y^h-y^{h-1}}^2+\frac{C}{2}\norm{x^{k+1}-y^k}^2\\
&= \left(\frac{L^2\tau}{2C} \sum_{h=k-\tau+1}^k(h-k+\tau)\norm{y^h-y^{h-1}}^2-\frac{L^2\tau}{2C} \sum_{h=k+1-\tau+1}^{k+1}(h-(k+1)+\tau)\norm{y^h-y^{h-1}}^2\right)\\
&~~~~~~~+\frac{L^2\tau^2}{2C}\norm{y^{k+1}-y^{k}}^2+\frac{C}{2}\norm{x^{k+1}-y^k}^2
\end{split}
\end{equation}
where $a)$ is due the triangle inequality and $d_k\leq \tau$. 
The linear extrapolation step for the
momentum acceleration in Algorithm \ref{AAPCD-Algo} yields  
\begin{equation}\label{deter-Ex-Acceleration-InE}
\begin{split}
\norm{y^{h+1}-y^{h}}&= (1+\beta_h) \E\norm{x^{h+1}-y^{h}}\qquad \text{for } h\in {\Gamma^c}_{k}^{r}, \beta_h>0\\
\norm{y^{h+1}-y^{h}}&= \E\norm{x^{h+1}-y^{h}}\qquad \text{for } h\in {\Gamma^0}_{k}^{r}\\
\norm{y^{h+1}-y^{h}}&= (1+\beta_h) \E\norm{x^{h+1}-y^{h}}\qquad \text{for } h\in {\Gamma}_{k}^{r}, \beta_h<0.
\end{split} 
\end{equation}
Thus, by taking total expectation on both sides of \eqref{stc-unbd-deay-lemm1-eq1} we have
\begin{equation}
\begin{split}
\E[\xi_{k+1}]&\leq\E[\xi_k]-\E\Iprod{{x}_{j_k}^{k+1}-y_{j_k}^{k}}{\nabla_{j_k} f(y^{k})-\nabla_{j_k} f(\hat{y}^{k})}\\
&~~~~~~~+\frac{L^2\tau^2}{2C}\E\norm{y^{k+1}-y^{k}}^2+\frac{C}{2}\E\norm{x^{k+1}-y^k}^2\\
&\leq\E[\xi_k]-\E\Iprod{{x}_{j_k}^{k+1}-y_{j_k}^{k}}{\nabla_{j_k} f(y^{k})-\nabla_{j_k} f(\hat{y}^{k})}\\
&~~~~~~~+\left(\frac{L^2\tau^2(1+\beta_k)^2}{2C}+\frac{C}{2}\right)\E\norm{x^{k+1}-y^k}^2.\\
\end{split}
\end{equation}
Therefore, we have
\begin{equation}
\begin{split}
\E[F(x^{k+1})+\xi_{k+1}]+&\frac{1}{2\eta}\E\norm{x^{k+1}-y^{k}}^2\\
\leq \E[F(y^{k})+\xi_{k}]+&\frac{L}{2}\E\norm{x^{k+1}-y^{k}}^2+\left(\frac{L^2\tau^2(1+\beta_k)^2}{2C}+\frac{C}{2}\right)\E\norm{x^{k+1}-y^k}^2.\\
\end{split}
\end{equation}
Hence, we can derive
\begin{equation}
\begin{split}
\E[F(x^{k+1})+\xi_{k+1}&]\\
\leq \E[F(y^{k})+&\xi_{k}]-\left(\frac{1}{2\eta}-\frac{L}{2}-\frac{L^2\tau^2(1+\beta_k)^2}{2C}-\frac{C}{2}\right)\E\norm{x^{k+1}-y^k}^2.\\
\end{split}
\end{equation}
By choosing $C={L\tau(1+\beta)}$, we obtain
\begin{equation}
\begin{split}
\E[F(x^{k+1})+\xi_{k+1}&]\\
\leq \E[F(y^{k})+&\xi_{k}]-\left(\frac{1}{2\eta}-\frac{L}{2}-{L\tau(1+\beta_k)}\right)\E\norm{x^{k+1}-y^k}^2\\
\end{split}
\end{equation}
and the result follows from the definition of $G(x^{k+1})$.
\end{proof}
\noindent\subsection{Proof of Theorem \ref{unbdd-P2-theo1}}
\begin{proof}
Applying Lemma \ref{unbdd-P2-lemma1} with $x=x^k$, $y=y^k$, we obtain that 
\begin{equation}\label{unbdd-P2-Theo1-eq1}
\begin{split}
G(x^{k+1})\leq G(y^{k})-\left(\frac{1}{2\eta_k}-\frac{L}{2}-\sqrt{{\delta_0}{\mu_{d_k}}}L(1+\beta_k)\right)\norm{x^{k+1}_{j_k}-y^k_{j_k}}^2.\\
\end{split}
\end{equation}
Since $\eta=\frac{c}{L+2L\sqrt{{\delta_0}{\mu_{T_1}}}(1+\beta)}$ and $-1 < \beta_k < {\frac{\sqrt{\mu_{T_1}}}{c\sqrt{\mu_{d_k}}}}(1+\beta)-1$, it follows that $G (x^{k+1})\leq G(y^k)$. Moreover, the update rule of the deterministic AAPCD guarantees that $F (y^{k+1})\leq F(x^{k+1})$ and hence we have $G (y^{k+1})\leq G(x^{k+1})$. In
summary, for all $k$ the following inequality holds:
\begin{equation}\label{unbdd-P2-Theo1-eq2}
G(y^{k+1})\leq G(x^{k+1})\leq G(y^{k})\leq G(x^{k}).
\end{equation}
Hence from \eqref{unbdd-P2-Theo1-eq1} we have
\begin{equation}\label{unbdd-P2-Theo1-eq11}
\begin{split}
G(y^{k+1})\leq G(y^{k})-\left(\frac{1}{2\eta_k}-\frac{L}{2}-\sqrt{{\delta_0}{\mu_{d_k}}}L(1+\beta_k)\right)\norm{x^{k+1}_{j_k}-y^k_{j_k}}^2.\\
\end{split}
\end{equation}
This equation shows $\norm{x^{h}-y^{h-1}}^2$ is summable. Thus, we have $\lim_{k} \norm{x^{k}-y^{k-1}}^2 = 0$.
Since $\sum_{h=1}^\infty \delta_{h} < \infty$, the series $A_k = \sum_{h=1}^k \delta_{k-h} \norm{y^{h}-y^{h-1}}^2$ converges to zero, i.e., $\xi_k\to 0$. Combing further with the fact that $F(x^k), F(y^k) \geq \inf F > -\infty$ for all $k$, we conclude that $\{G(x^k )\}, \{G(y^k)\}$ converge to the same limit $F^*$, i.e.,
\begin{equation}
\lim_{k\to\infty}F(x^k) = \lim_{k\to\infty}F(y^k) =\lim_{k\to\infty}G(x^k) = \lim_{k\to\infty}G(y^k)= F^*.
\end{equation}
On the other hand, by induction we conclude from equation \eqref{unbdd-P2-Theo1-eq2} that for all $k$
\[
F(y^k)\leq G(y^k)\leq F(x^0)\qquad F(x^k)\leq G(x^k)\leq F(x^0).
\]
Combining with Assumption \ref{P2-assu1} that $F$ has bounded sublevel set, we conclude that $\{x^k\}$ and $\{y^k\}$ are bounded and thus
have bounded limit points.

\noindent Since $\norm{y^k - x^{k+1}} \to 0$, ${x^k}$ and ${y^k}$ share the same set of limit points $\Omega$ which is compact in $\mathbb{R}^m$. We fix any limit point $x\in\Omega$, say $x^{k_t+1} \to x$, $y^{k_t} \to x$. Note that the continuity of $f$ yields $\lim_{t\to\infty} f(x^{k_t}) = f(x)$. 
By the definition of $x^{k+1}_j$ as a proximal point, we have 
\begin{equation}
\begin{split}
&\Iprod{\nabla_j f(\hat{y}^{l(k_t+1,j)})}{x^{l(k_t+1,j)+1}_j-y^{l(k_t+1,j)}_j}+\frac{1}{2\eta}\norm{x^{l(k_t+1,j)+1}_j-y^{l(k_t+1,j)}_j}^2+g_j(x^{l(k_t+1,j)+1}_j)\\
&\leq\Iprod{\nabla_j f(\hat{y}^{l(k_t+1,j)})}{x_j-y^{l(k_t+1,j)}_j} + \frac{1}{2\eta}\norm{x_j-y^{l(k_t+1,j)}_j}^2+g_j(x_j).
\end{split}
\end{equation}
 Taking $\lim\sup$ on both sides and note that $x^{k_t+1} - y^{k_t} \to 0$, $y^{k_t} \to x$, we obtain that $\lim\sup_{t\to\infty} g_j(x^{k_t+1}_j) \leq g_j(x_j)$. Since
$g$ is lower semicontinuous and $x^{k_t}\to x$, it follows that $\lim\sup_{t\to\infty} g_j(x^{k_t+1}_j)\geq g_j(x_j)$. 
Combining both inequalities, we
conclude that $\lim_{t\to\infty} g_j(x^{k_t+1}_j) = g_j(x_j)$. We then obtain
$\lim_{t\to\infty} F(x^{k_t}) = F(x)$. Since $\lim_{t\to\infty} F(x^{k_t}) = F^*$ by equation \eqref{unbdd-P2-Theo1-eq2}, we have
\begin{equation}
F(x)=F^*,\qquad \forall x\in\Omega.
\end{equation}
Thus, $F$ remains constant on the set of limit points $\Omega$. 
By optimality condition of the proximal gradient step of AAPCD, we obtain that
\begin{equation}\label{unbdd-P2-Theo1-eq4}
\begin{split}
-\nabla_{j_k} f(\hat{y}^{k})-\frac{1}{\eta_k}(x_{j_k}^{k+1}-y_{j_k}^k)&\in \partial g_{j_k}(x^{k+1}_{j_k})\\
\leftrightarrow\underbrace{\nabla_{j_k} f(x^{k+1})-\nabla_{j_k} f(\hat{y}^{k})-\frac{1}{\eta_k}(x^{k+1}_{j_k}-y^k_{j_k})}_{q^{k+1}_{j_k}}&\in \partial_{x_{j_k}} F(x^{k+1})
\end{split}
\end{equation}
and for $j\neq j_k$,
\begin{equation}\label{unbdd-P2-Theo1-eq6}
\begin{split}
-\nabla_{j} f(\hat{y}^{l(k,j)})-\frac{1}{\eta_k}(x_{j}^{l(k,j)+1}-y_{j}^{l(k,j)})&\in \partial g_{j}(x^{k+1}_{j})\\
\leftrightarrow\underbrace{\nabla_{j} f(x^{k+1})-\nabla_{j} f(\hat{y}^{l(k,j)})-\frac{1}{\eta_k}(x^{l(k,j)+1}_{j}-y^{l(k,j)}_{j})}_{q^{k+1}_{j}}&\in \partial_{x_{j}} F(x^{k+1}).
\end{split}
\end{equation}
We have also 
\begin{equation}\label{unbdd-P2-Theo1-eq61}
\frac{\partial y^{k+1}_{j_k}}{\partial x_{j_k}^{k+1}}\frac{\partial G(x^{k+1})}{\partial y^{k+1}_{j_k}} = \frac{\partial y^{k+1}_{j_k}}{\partial x_{j_k}^{k+1}}\sqrt{\frac{\mu_{d_k}}{\delta_0}}\frac{L}{(1+\beta_k)} \delta_0 (y^{k+1}_{j_k}-y^{k}_{j_k})
\end{equation}
where $\norm{\frac{\partial y^{k+1}}{\partial x^{k+1}}} \leq (1+\beta_k)$. For $k+1\in S_T$ we have
\begin{equation}\label{unbdd-P2-Theo1-eq8}
\begin{split}
\norm{q_j^{k+1}}^2 &\myeq{a)}{\leq} \left(\norm{\nabla_j f(x^{k+1})-\nabla_j f(\hat{y}^{l(k+1,j)})-\frac{1}{\eta}(x^{k+1}_j-y^{l(k+1,j)}_j)}\right)^2\\
&\myeq{b)}{\leq} \left(\frac{1}{\eta} \norm{x^{k+1}_j-y^{{l(k+1,j)}}_j}+L\norm{x^{k+1}_j-\hat{y}^{l(k+1,j)}_j}\right)^2\\
&\myeq{c)}{\leq} \left((\frac{1}{\eta}+L)\norm{x^{{k+1}}_j-y^{{l(k+1,j)}}_j}+L\norm{y^{l(k+1,j)}_j-\hat{y}^{l(k+1,j)}_j}\right)^2\\
&\leq 2(\frac{1}{\eta}+L)^2\norm{x^{{k+1}}_j-y^{{l(k+1,j)}}_j}^2+2L^2\norm{y^{l(k+1,j)}_j-\hat{y}^{l(k+1,j)}_j}^2\\
&\myeq{d)}{\leq} 2(\frac{1}{\eta}+L)^2\norm{x^{{l(k+1,j)+1}}_j-y^{{l(k+1,j)}}_j}^2+2L^2T\sum_{h=l(k+1,j)-T+1}^{l(k+1,j)}\norm{y^h_j-y^{h-1}_j}^2\\
&\myeq{e)}{\leq} 2(\frac{1}{\eta}+L)^2K\sum_{h=k+1-K}^{{k+1}}\norm{x^h-y^{h-1}}^2+2L^2T\sum_{h=k+1-T-K}^{{k+1}}\norm{y^h-y^{h-1}}^2
\end{split}
\end{equation}
where $a)$ follows from \eqref{unbdd-P2-Theo1-eq4} and \eqref{unbdd-P2-Theo1-eq6}, $b)$ by the Lipschitz of $\nabla f$,  
$c)$ is by the triangle inequality, $d)$ from triangular inequality and  and the fact that $k+1\in S_T$ and $e)$ is due to Assumption \ref{unbdd-determinstic-rule}.
Thus, from this we have $\norm{q_j^{k+1}}\to 0$.
We have shown $x^{k_t}\to x$, $F (x^{k_t})\to F(x)$ and that
$\partial F(x^{k_t})$ converges to $0$. Therefore, $0\in\partial F(x)$ for all $x\in\Omega$. 
\end{proof}
\noindent\subsection{Proof of Theorem \ref{unbdd-deter-subseq-conver}}
\begin{proof}
 From equation \eqref{unbdd-P2-Theo1-eq11} we have
\begin{equation}
\begin{split}
G(y^{h+1})\leq G(y^{h})-\left(\frac{1}{2\eta_h}-\frac{L}{2}-\sqrt{{\delta_0}{\mu_{d_h}}}L(1+\beta_h)\right)\norm{x^{h+1}-y^h}^2.\\
\end{split}
\end{equation}
Recall $\eta=\frac{c}{L+2L\sqrt{{\delta_0}{\mu_{T_1}}}(1+\beta)}$ for $c\in (0,1)$, and $-1 < \beta_k < {\frac{\sqrt{\mu_{T_1}}}{c\sqrt{\mu_{d_k}}}}(1+\beta)-1$ implies
\begin{equation}
\begin{split}
(\frac{1}{c}-1)&\frac{L}{2}\norm{x^{h+1}-y^{h}}^2 \leq G(y^{h})-G(y^{h+1}).\\
\end{split}
\end{equation}
By summing the above inequality over ${k+1-T-K},\ldots, k$ iterations we obtain
\begin{equation}\label{unbdd-P2-Theo2-eq1}
\begin{split}
(\frac{1}{c}-1)&\frac{L}{2}\sum_{h=k+1-T-K}^{k}\norm{x^{h+1}-y^{h}}^2 \leq G(y^{k+1-T-K})-G(y^{k+1}).\\
\end{split}
\end{equation}
We have shown in Theorem \ref{unbdd-P2-theo1} that $F(x^k)\to F^*$, and it is also clear that $\text{dist}_{\Omega}(x^k)\to 0$. Thus, for any $\epsilon,\gamma > 0$ there is $k_0$ such that for all
$k \geq k_0$, we have
\begin{equation}
x^k\in\{x|\,\,\text{dist}_{\Omega}(x)\leq \epsilon, F^* < F(x) < F^*+\gamma\}.
\end{equation}
Since $\Omega$ is compact and $F$ is constant on it, the Uniformized KL property implies that for all $k \geq k_0$ 
\begin{equation}\label{deter-P2-Theo2-eq3}
\phi'(F(x^{k+1})-F^*)\text{dist}_{\partial F(x^{k+1})}(0)\geq 1.
\end{equation}
Recall that ${r_k := F(y^k) - F^*}$. Then, we have
\begin{equation}
\begin{split}
1&\myeq{a)}{\leq} \left(\phi'(F(x^{k+1})-F^*)\text{dist}_{\partial F(x^{k+1})}(0)\right)^2\\
 &\myeq{b)}{\leq} \left(\phi'(r_{k+1})\text{dist}_{\partial F(x^{k+1})}(0)\right)^2\\
 &\myeq{c)}{\leq} (\phi'(r_{k+1}))^2\bigg(2(\frac{1}{\eta}+L)^2\sum_{j=1}^m\norm{x^{{l(k+1,j)+1}}_j-y^{{l(k+1,j)}}_j}^2\\&~~~+2L^2T_1\sum_{j=1}^m\sum_{h=l(k+1,j)-T_1+1}^{l(k+1,j)}\norm{y^h_j-y^{h-1}_j}^2\\
&~~~+2L^2T\sum_{j=1}^m\sum_{h=l(k+1,j)-T+1}^{l(k+1,j)}\norm{y^h_j-y^{h-1}_j}^2\bigg)\\
&{\leq} (\phi'(r_{k+1}))^2\bigg( 2(\frac{1}{\eta}+L)^2\sum_{h=k+1-K}^{{k+1}}\norm{x^h-y^{h-1}}^2+2L^2T_1\sum_{h=k+1-T_1-K}^{{k+1}}\norm{y^h-y^{h-1}}^2\\
&~~~+2L^2T\sum_{h=k+1-T-K}^{{k+1}}\norm{y^h-y^{h-1}}^2\bigg)\\
&{\leq} (\phi'(r_{k+1}))^2\bigg( 2(\frac{1}{\eta}+L)^2\sum_{h=k+1-K}^{{k+1}}\norm{x^h-y^{h-1}}^2+{2L^2T_1}(1+\beta)^2\sum_{\substack{ h\in {\Gamma^c}_{k+1-{K}}^{k+1}}}\norm{x^h-y^{h-1}}^2\\
&+{2L^2T}(1+\beta'')^2\sum_{\substack{ h\in {\Gamma^0}_{k+1-T-{K}}^{k+1}\cup{\Gamma}_{k+1-T-{K}}^{k+1}}}\norm{x^h-y^{h-1}}^2\bigg)\\
&\myeq{d)}{\leq} (\phi'(r_{k+1}))^2\left(\frac{2(\frac{1}{\eta}+L)^2+3(1+\beta)^2L^2T_1+2(1+\beta'')^2L^2T}{(\frac{1}{c}-1)\frac{L}{2}}\right)\\
&~~~~\times(G(y^{k+1-T-K})-G(y^{k+1}))\\
&\myeq{e)}{\leq} (\phi'(r_{k+1}))^2b_1(G(y^{k+1-T-K})-G(y^{k+1}))\\
&\myeq{f)}{\leq} (\phi'(r_{k+1}))^2b_1(F(y^{0})-F(y^{k+1}))
\end{split}
\end{equation}
where $a)$ follows from \eqref{deter-P2-Theo2-eq3}, $b)$ is due to $r_{k+1}\leq F(x^{k+1})-F^*$ and the fact that $\phi'$ is nonincreasing, $c)$ from \eqref{unbdd-P2-Theo1-eq8}, $d)$ is a direct
computation using \eqref{unbdd-P2-Theo2-eq1} and $\beta'' = \max\{\beta',0\}$, $e)$ is also a result of the definition of $b_1$ and $f)$ by $G(y^{k+1-T-K})\leq F(y^0)$ and $F(y^{k+1})\leq G(y^{k+1})$.
 We have that $\phi'(t)={e}t^{\theta-1}$. Thus the above inequality becomes 
\begin{equation}\label{deter-p2-rate-relation-LK}
1\leq b_1e^2 r_{k+1}^{2\theta-2}(r_{{0}}-r_{k+1})\qquad k+1\in S_T.
\end{equation}
\noindent {\bf Part} $1:$ Suppose that $\theta = 1$, then for all $k$, we have $r_{0}-r_{k+1}\geq \frac{1}{b_1e^2}>0$, which cannot hold because $r_0 < \frac{1}{b_1e^2}$. Thus, $\{F(y^k)\}$ must converge in finitely many steps, which is by Theorem \ref{unbdd-P2-theo1} the stationary point of $F$.

In the following  we assume that $r_k\neq 0$ for all $k\in S_T$ because otherwise the algorithm terminates.

\noindent {\bf Part} $2:$ Suppose that $\theta\in[\frac{1}{2},1)$. We select $k_1$ large enough such that $r_{k+1}^{2-2\theta} \geq r_{k+1}$, for all $k\geq k_1$. Then for all $k \geq k_1$, and $k+1\in S_T$ we have  
\begin{equation}\label{undbb-deter-part2-eq1}
r_{k+1}\leq \frac{b_1e^2}{1+b_1e^2} r_{0}.
\end{equation}

\noindent {\bf Part} $3:$ Suppose that $\theta\in (0,\frac{1}{2})$. Let $h(s):=s^{2\theta-2}$. Then from \eqref{deter-p2-rate-relation-LK} we find that 
\begin{equation}
\begin{split}
\frac{1}{b_1e^2}&\leq h(r_{k+1})(r_{0}-r_{k+1}) = \frac{h(r_{k+1})}{h(r_{0})} h(r_{0})(r_{0}-r_{k+1})\\
&\leq \frac{h(r_{k+1})}{h(r_{0})}\int_{r_{k+1}}^{r_{0}} h(s) ds = \frac{h(r_{k+1})}{h(r_{0})}\frac{r_{k+1}^{2\theta-1}-r_{0}^{2\theta-1}}{1-2\theta}.
\end{split}
\end{equation}
Let $R\in(1, \infty)$ be a fixed number. We consider two cases:

\noindent {\bf Case 1:} Let $\frac{h(r_{k+1})}{h(r_{0})}\leq R$. 
Then we have 
\begin{equation}\label{unbdd-KL-case3case1}
\begin{split}
\frac{1}{b_1e^2R}&\leq \frac{r_{k+1}^{2\theta-1}-r_{{0}}^{2\theta-1}}{1-2\theta}.
\end{split}
\end{equation}

\noindent {\bf Case 2:} Let $\frac{h(r_{k+1})}{h(r_{0})}\geq R$. 
Then, since $h(r_{k+1})\geq h(r_{0}) R$,  we have $r_{k+1}^{2\theta-1}\geq r_{0}^{2\theta-1} R^{\frac{2\theta-1}{2\theta-2}}$. Therefore, we obtain
\[
 r_{k+1}^{2\theta-1}\geq r_{0}^{2\theta-1} R^{\frac{2\theta-1}{2\theta-2}}.
\]
Thus, we can deduce that
\begin{equation}\label{unbdd-KL-case3case2}
\begin{split}
\frac{r_{0}^{2\theta-1} (R^{\frac{2\theta-1}{2\theta-2}}-1)}{1-2\theta}\leq \frac{r_{k+1}^{2\theta-1}-r_{0}^{2\theta-1}}{1-2\theta}.
\end{split}
\end{equation}
Combining equations \eqref{unbdd-KL-case3case1} and \eqref{unbdd-KL-case3case2} yields 
\begin{equation}\label{undbb-deter-part3-eq1}
\begin{split}
\min(\frac{1}{b_1e^2R},\frac{r_{0}^{2\theta-1} (R^{\frac{2\theta-1}{2\theta-2}}-1)}{1-2\theta}) \leq \frac{r_{k+1}^{2\theta-1}-r_{0}^{2\theta-1}}{1-2\theta}.
\end{split}
\end{equation}
Hence we have
\begin{equation}
\begin{split}
b_2 \leq \frac{r_{k+1}^{2\theta-1}-r_{0}^{2\theta-1}}{1-2\theta},
\end{split}
\end{equation}
where $b_2 = \min(\frac{1}{b_1e^2R},\frac{r_{0}^{2\theta-1} (R^{\frac{2\theta-1}{2\theta-2}}-1)}{1-2\theta})$. Thus, we have the following bound for all $k+1\in S_T$
\begin{equation}
r_{k+1}\leq \left(\frac{1}{b_2{({1-2\theta})}+r_{{0}}^{2\theta-1}}\right)^{\frac{1}{1-2\theta}}.
\end{equation}
\end{proof}
\noindent\subsection{Proof of Theorem \ref{bd-deter-corro}}
\begin{proof}
We have shown in Theorem \ref{unbdd-P2-theo1} that $\tilde{G}(x^k)\to F^*$, and it is also clear that $\text{dist}_{\Omega}(x^k)\to 0$. Furthermore, similar to the proof of Theorem \ref{unbdd-P2-theo1}, we can show the elements of $\Omega$ are the critical points of $\tilde{G}$. Thus, for any $\epsilon,\gamma > 0$ there is $k_0$ such that for all
$k \geq k_0$, we have
\begin{equation}
x^k\in\{x|\,\,\text{dist}_{\Omega}(x)\leq \epsilon, F^* < \tilde{G}(x) < F^*+\gamma\}.
\end{equation}
Since $\Omega$ is compact and $\tilde{G}$ is constant on it, the Uniformized KL property implies that for all $k \geq k_0$ 
\begin{equation}\label{bdd-deter-P2-Theo2-eq3}
\phi'(\tilde{G}(x^{k+1})-F^*)\text{dist}_{\partial \tilde{G}(x^{k+1})}(0)\geq 1.
\end{equation}
Recall that ${r_k := \tilde{G}(y^k) - F^*}$. Then, we have
\begin{equation}
\begin{split}
1&\myeq{a)}{\leq} \left(\phi'(\tilde{G}(x^{k+1})-F^*)\text{dist}_{\partial \tilde{G}(x^{k+1})}(0)\right)^2\\
 &\myeq{b)}{\leq} \left(\phi'(r_{k+1})\text{dist}_{\partial \tilde{G}(x^{k+1})}(0)\right)^2\\
 &\myeq{c)}{\leq} (\phi'(r_{k+1}))^2\bigg(3(\frac{1}{\eta}+L)^2\sum_{j=1}^m\norm{x^{{l(k+1,j)+1}}_j-y^{{l(k+1,j)}}_j}^2\\&~~~+3L^2T_1\sum_{j=1}^m\sum_{h=l(k+1,j)-T_1+1}^{l(k+1,j)}\norm{y^h_j-y^{h-1}_j}^2\\
&~~~+3L^2T\sum_{j=1}^m\sum_{h=l(k+1,j)-T+1}^{l(k+1,j)}\norm{y^h_j-y^{h-1}_j}^2+{6L^2\mu_{T}} \delta_{0}\norm{y^{k+1}-y^{k}}^2\bigg)\\
&{\leq} (\phi'(r_{k+1}))^2\bigg( 3(\frac{1}{\eta}+L)^2\sum_{h=k+1-K}^{{k+1}}\norm{x^h-y^{h-1}}^2+3L^2T_1\sum_{h=k+1-T_1-K}^{{k+1}}\norm{y^h-y^{h-1}}^2\\
&~~~+3L^2T\sum_{h=k+1-T-K}^{{k+1}}\norm{y^h-y^{h-1}}^2+{6L^2\mu_{T}} \delta_{0}\norm{y^{k+1}-y^{k}}^2\bigg)\\
&{\leq} (\phi'(r_{k+1}))^2\bigg( 3(\frac{1}{\eta}+L)^2\sum_{h=k+1-K}^{{k+1}}\norm{x^h-y^{h-1}}^2+{3L^2T_1}(1+\beta)^2\sum_{\substack{ h\in {\Gamma^c}_{k+1-{K}}^{k+1}}}\norm{x^h-y^{h-1}}^2\\
&+{3L^2T}(1+\beta'')^2\sum_{\substack{ h\in {\Gamma^0}_{k+1-T-{K}}^{k+1}\cup{\Gamma}_{k+1-T-{K}}^{k+1}}}\norm{x^h-y^{h-1}}^2+{6L^2\mu_{T}} \delta_{0}(1+\beta)^2\norm{x^{k+1}-y^{k}}^2\bigg)\\
&\myeq{d)}{\leq} (\phi'(r_{k+1}))^2\left(\frac{3(\frac{1}{\eta}+L)^2+3(1+\beta)^2L^2T_1+3(1+\beta'')^2L^2T+6(1+\beta)^2{L^2\mu_{T}} \delta_{0}}{(\frac{1}{c}-1)\frac{L}{2}}\right)\\
&~~~~\times(\tilde{G}(y^{k+1-T-K})-\tilde{G}(y^{k+1}))\\
&\myeq{e)}{\leq} (\phi'(r_{k+1}))^2b_1(\tilde{G}(y^{k+1-T-K})-\tilde{G}(y^{k+1}))
\end{split}
\end{equation}
where $a)$ follows from \eqref{bdd-deter-P2-Theo2-eq3}, $b)$ is due to $r_{k+1}\leq \tilde{G}(x^{k+1})-F^*$ and the fact that $\phi'$ is nonincreasing, $c)$ from \eqref{unbdd-P2-Theo1-eq8} and \eqref{unbdd-P2-Theo1-eq61}, $d)$ is a direct
computation using \eqref{unbdd-P2-Theo2-eq1} and $\beta'' = \max\{\beta',0\}$, and  $e)$ is also a result of the definition of $b_1$.
 We have that $\phi'(t)={e}t^{\theta-1}$. 
Thus, we obtain 
\begin{equation}\label{bdd-deter-p2-rate-relation-LK}
1\leq b_1e^2 r_{k+1}^{2\theta-2}(r_{{k+1-T-K}}-r_{k+1}).
\end{equation}

\noindent {\bf Part} $1:$ Suppose that $\theta = 1$, then for all $k$, by \eqref{bdd-deter-p2-rate-relation-LK} we have $r_{k+1-T-K}-r_{k+1}\geq \frac{1}{b_1e^2}>0$, which cannot hold because $r_{k}\to 0$. Thus, $\{G(y^k)\}$ must converge in finitely many steps, which is the stationary point of $\tilde{G}$.

In the following  we assume that $r_k\neq 0$ for all $k\in S_T$ because otherwise the algorithm terminates.

\noindent {\bf Part} $2:$ Suppose that $\theta\in[\frac{1}{2},1)$. We select $k_1$ large enough such that $r_{k+1}^{2-2\theta} \geq r_{k+1}$, for all $k\geq k_1$. Then for all $k \geq k_1$, and $k+1\in S_T$ we have  
\begin{equation}
r_{k+1}\leq \frac{b_1e^2}{1+b_1e^2} r_{k+1-T-K}\leq \left(\frac{b_1e^2}{1+b_1e^2}\right)^{\lfloor\frac{k+1-k_1}{T+K}\rfloor}r_{k_1}.
\end{equation}

\noindent {\bf Part} $3:$ Suppose that $\theta\in (0,\frac{1}{2})$. Let $h(s):=s^{2\theta-2}$. Then from \eqref{bdd-deter-p2-rate-relation-LK} we find that 
\begin{equation}
\begin{split}
\frac{1}{b_1e^2}&\leq h(r_{k+1})(r_{k+1-T-K}-r_{k+1}) = \frac{h(r_{k+1})}{h(r_{k+1-T-K})} h(r_{k+1-T-K})(r_{k+1-T-K}-r_{k+1})\\
&\leq \frac{h(r_{k+1})}{h(r_{k+1-T-K})}\int_{r_{k+1}}^{r_{k+1-T-K}} h(s) ds = \frac{h(r_{k+1})}{h(r_{k+1-T-K})}\frac{r_{k+1}^{2\theta-1}-r_{k+1-T-K}^{2\theta-1}}{1-2\theta}.
\end{split}
\end{equation}
Let $R\in(1, \infty)$ be a fixed number. We consider two cases:

\noindent {\bf Case 1:} Let $\frac{h(r_{k+1})}{h(r_{k+1-T-K})}\leq R$. 
Then we have 
\begin{equation}\label{bdd-unbdd-KL-case3case1}
\begin{split}
\frac{1}{b_1e^2R}&\leq \frac{r_{k+1}^{2\theta-1}-r_{{k+1-T-K}}^{2\theta-1}}{1-2\theta}.
\end{split}
\end{equation}

\noindent {\bf Case 2:} Let $\frac{h(r_{k+1})}{h(r_{k+1-T-K})}\geq R$. 
Then, since $h(r_{k+1})\geq h(r_{k+1-T-K}) R$,  we have $r_{k+1}^{2\theta-1}\geq r_{k+1-T-K}^{2\theta-1} R^{\frac{2\theta-1}{2\theta-2}}$. Therefore, we obtain
\[
 r_{k+1}^{2\theta-1}\geq r_{k+1-T-K}^{2\theta-1} R^{\frac{2\theta-1}{2\theta-2}}.
\]
Thus, we can deduce that
\begin{equation}\label{bdd-unbdd-KL-case3case2}
\begin{split}
\frac{r_{k+1-T-K}^{2\theta-1} (R^{\frac{2\theta-1}{2\theta-2}}-1)}{1-2\theta}\leq \frac{r_{k+1}^{2\theta-1}-r_{k+1-T-K}^{2\theta-1}}{1-2\theta}.
\end{split}
\end{equation}
Combining equations \eqref{bdd-unbdd-KL-case3case1} and \eqref{bdd-unbdd-KL-case3case2} yields 
\begin{equation}\label{bdd-undbb-deter-part3-eq1}
\begin{split}
\min(\frac{1}{b_1e^2R},\frac{r_{k+1-T-K}^{2\theta-1} (R^{\frac{2\theta-1}{2\theta-2}}-1)}{1-2\theta}) \leq \frac{r_{k+1}^{2\theta-1}-r_{k+1-T-K}^{2\theta-1}}{1-2\theta}.
\end{split}
\end{equation}
From equation \eqref{bdd-undbb-deter-part3-eq1} we have
\begin{equation}
\begin{split}
tb_2 \leq \frac{r_{k+1+t{(T+K)}}^{2\theta-1}-r_{k+1}^{2\theta-1}}{1-2\theta},
\end{split}
\end{equation}
where $b_2 = \min(\frac{1}{b_1e^2R},\frac{r_{0}^{2\theta-1} (R^{\frac{2\theta-1}{2\theta-2}}-1)}{1-2\theta})$. Thus, we have the following bound for all $k \geq k_0$
\begin{equation}
r_{k+1}\leq \left(\frac{1}{{\lfloor\frac{k+1-k_0}{T+K}\rfloor}b_2{({1-2\theta})}+r_{{0}}^{2\theta-1}}\right)^{\frac{1}{1-2\theta}}.
\end{equation}
\end{proof}
\noindent\subsection{Additional Numerical Results}
\begin{figure*}[htbp]
\centering
\subfloat[a]{
\centering
\includegraphics[width=0.40\linewidth]{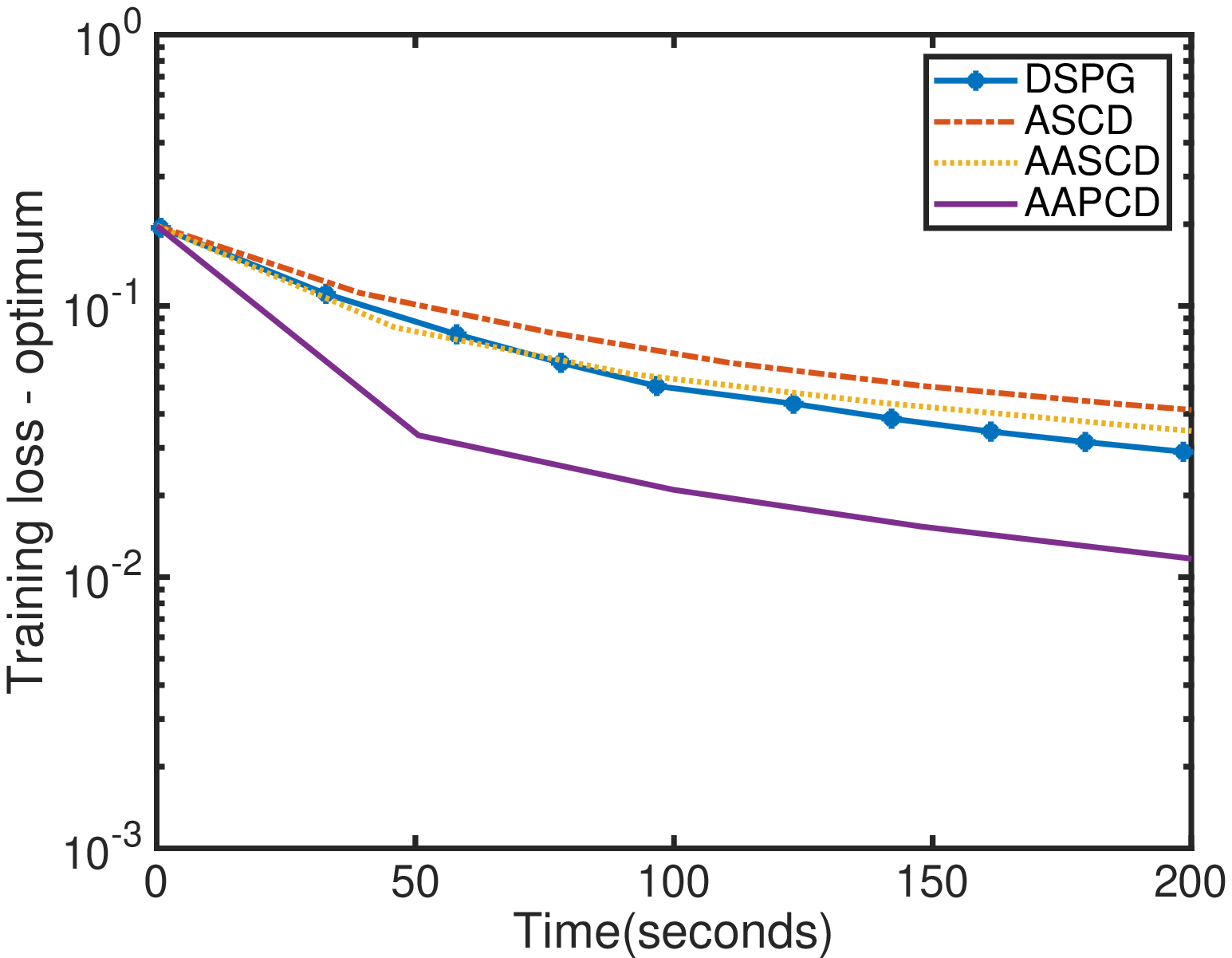}}%
\hfill
\subfloat[b]{
\centering
\includegraphics[width=0.40\linewidth]{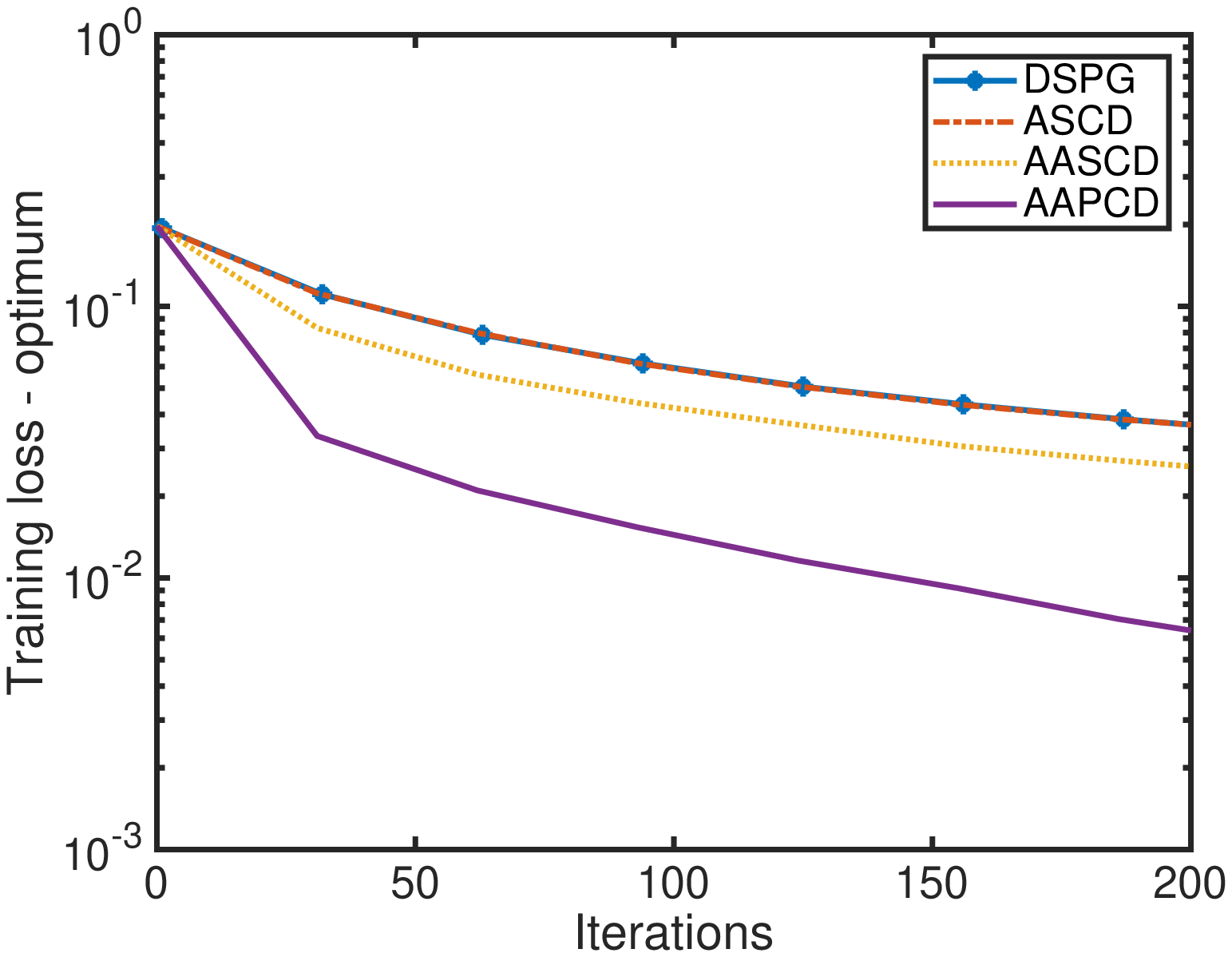}}%
\setlength{\abovecaptionskip}{2pt}
\caption{Training loss residual  versus time (Figure(a)) and iteration (Figure(b)) plot of AAPCD, ASCD, AASCD, and DSGD. }
\label{fig:algo_comp_supp}
\end{figure*}

We focus here on Sigmoid loss with $l_2$ regularization term:
\[
f(x) = \frac{1}{n}\sum_{i=1}^n \frac{1}{1+e^{b_ix^Ta_i}}
\]
and
\[
g(x) = \lambda \sum_{j=1}^m\norm{x_j}.
\]
Experiments are performed on $\it{Covetype}$ dataset. We conduct experiments for comparing AAPCD with other algorithms: ASCD, AASCD, and DSPG. For all the experiments we set the number of workers P=32. The block size for all experiments is $10$. The convergence results are presented in Figure \ref{fig:algo_comp_supp}. The results from AAPCD outperforms other algorithms. Moreover AAPCD obtains a much smaller objective value. Note that DSPG and ASCD have similar performance. AASCD is faster than ASCD, but the analyses for AASCD is only developed for convex functions. Overall, the results show that our algorithm is very efficient for nonconvex functions such as Sigmoid loss.
